\documentclass[a4paper,leqno,final]{amsart}
\usepackage[foot]{amsaddr}
\RequirePackage{ifthen}
\newboolean{omarfont}
\setboolean{omarfont}{false}
\newboolean{usesvjour}
\setboolean{usesvjour}{false}%
    \RequirePackage{ifthen}
    \provideboolean{useutopia}
    \setboolean{useutopia}{false}%
    \provideboolean{usesvjour}
    \setboolean{usesvjour}{false}%
    \ifthenelse{\boolean{useutopia}}{%
      \RequirePackage{ifthen}
      \RequirePackage{iftex}
      \provideboolean{usemathrsfs}%
      \provideboolean{useutopia}
      \setboolean{useutopia}{false}%
      \ifXeTeX
        \RequirePackage[american,british]{babel}
        \RequirePackage{xltxtra}%
        \DeclareSymbolFont{usualmathcal}{OMS}{cmsy}{m}{n}
        \DeclareSymbolFontAlphabet{\mathcalbf}{usualmathcal}
        \ifthenelse{\boolean{useutopia}}{
          \RequirePackage[utopia,euro]{mathdesign}
          \RequirePackage[OMLmathbf,OMLmathsfit]{isomath}
        }{
          \RequirePackage{amssymb}    %
          \RequirePackage{bbold}    %
        }
        \RequirePackage{mathrsfs} %
        \setboolean{usemathrsfs}{true}%
      \else%
        \RequirePackage[utf8]{inputenc}%
        \ifthenelse{\boolean{useutopia}}{%
          \RequirePackage[utopia,euro]{mathdesign}%
          \RequirePackage[OMLmathrm,OMLmathbf,OMLmathsfit]{isomath}
          \DeclareSymbolFont{usualmathcal}{OMS}{cmsy}{m}{n}
          \DeclareSymbolFontAlphabet{\mathcalbf}{usualmathcal}
          \RequirePackage{stmaryrd} %
          \providecommand{\diracdelta}[1][]{\ensuremath{\deltaup_{#1}}}
          
          \providecommand{\lap}{\ensuremath{\Deltaup}}
          \providecommand{\measure}[1]{\ensuremath{\mathcalbf{\uppercase{#1}}}}
        }{%
          \RequirePackage{amssymb}    %
          \RequirePackage{mathrsfs} %
          \RequirePackage{bbold}    %
          \RequirePackage{stmaryrd} %
          \setboolean{usemathrsfs}{true}%
          \providecommand{\mathcalbf}{\mathcal}
        }%
        \RequirePackage[american,british]{babel}
      \fi%
      \RequirePackage{mathtools}%
      \RequirePackage{nicefrac}
      \RequirePackage{stmaryrd} %
      \RequirePackage{amsthm}%
      \RequirePackage[shortlabels]{enumitem}%
      \RequirePackage{xspace}%
      \RequirePackage{verbatim}%
      \RequirePackage{fancyvrb}
      \RequirePackage{listings}%
      \RequirePackage{xcolor}
      \RequirePackage{epsfig}
      \RequirePackage{graphicx}
      \RequirePackage{booktabs}%
      \RequirePackage{tikz}%
      \usetikzlibrary{calc}%
      \provideboolean{usebeamer}%
      \provideboolean{isthesis}%
      \provideboolean{isamsltex}%
      \provideboolean{issiamltex}%
    }{%
      \RequirePackage{ifthen}
      \RequirePackage{iftex}
      \provideboolean{usemathrsfs}%
      \provideboolean{useutopia}
      \setboolean{useutopia}{false}%
      \ifXeTeX
        \RequirePackage[american,british]{babel}
        \RequirePackage{xltxtra}%
        \DeclareSymbolFont{usualmathcal}{OMS}{cmsy}{m}{n}
        \DeclareSymbolFontAlphabet{\mathcalbf}{usualmathcal}
        \ifthenelse{\boolean{useutopia}}{
          \RequirePackage[utopia,euro]{mathdesign}
          \RequirePackage[OMLmathbf,OMLmathsfit]{isomath}
        }{
          \RequirePackage{amssymb}    %
          \RequirePackage{bbold}    %
        }
        \RequirePackage{mathrsfs} %
        \setboolean{usemathrsfs}{true}%
      \else%
        \RequirePackage[utf8]{inputenc}%
        \ifthenelse{\boolean{useutopia}}{%
          \RequirePackage[utopia,euro]{mathdesign}%
          \RequirePackage[OMLmathrm,OMLmathbf,OMLmathsfit]{isomath}
          \DeclareSymbolFont{usualmathcal}{OMS}{cmsy}{m}{n}
          \DeclareSymbolFontAlphabet{\mathcalbf}{usualmathcal}
          \RequirePackage{stmaryrd} %
          \providecommand{\diracdelta}[1][]{\ensuremath{\deltaup_{#1}}}
          
          \providecommand{\lap}{\ensuremath{\Deltaup}}
          \providecommand{\measure}[1]{\ensuremath{\mathcalbf{\uppercase{#1}}}}
        }{%
          \RequirePackage{amssymb}    %
          \RequirePackage{mathrsfs} %
          \RequirePackage{bbold}    %
          \RequirePackage{stmaryrd} %
          \setboolean{usemathrsfs}{true}%
          \providecommand{\mathcalbf}{\mathcal}
        }%
        \RequirePackage[american,british]{babel}
      \fi%
      \RequirePackage{mathtools}%
      \RequirePackage{nicefrac}
      \RequirePackage{stmaryrd} %
      \RequirePackage{amsthm}%
      \RequirePackage[shortlabels]{enumitem}%
      \RequirePackage{xspace}%
      \RequirePackage{verbatim}%
      \RequirePackage{fancyvrb}
      \RequirePackage{listings}%
      \RequirePackage{xcolor}
      \RequirePackage{epsfig}
      \RequirePackage{graphicx}
      \RequirePackage{booktabs}%
      \RequirePackage{tikz}%
      \usetikzlibrary{calc}%
      \provideboolean{usebeamer}%
      \provideboolean{isthesis}%
      \provideboolean{isamsltex}%
      \provideboolean{issiamltex}%
    }
    \colorlet{a}{magenta}
    \colorlet{b}{green!75!blue}
    \colorlet{c}{yellow!87.5!red}
    \colorlet{d}{cyan}
    \colorlet{e}{red}
    \colorlet{f}{blue}
    \colorlet{g}{white}
    \colorlet{h}{black!50}
    \colorlet{i}{black}
    \colorlet{j}{black!75}

    \RequirePackage{hyperref}
    \providecommand{\linkedurl}[1]{\url{#1}}
    \providecommand{\linkedemail}[1]{\href{mailto:#1}{#1}}
    
    \providecommand{\email}[1]{{\linkedemail{#1}}}

    \providecommand{\Ignore}[1]{}
    \providecommand{\ignore}[1]{}
    \providecommand{\freeze}[1]{}%
    \providecommand{\crossout}[1]{{\textcolor{i!20}{#1}}}
    \providecommand{\highlightcolor}{a}
    \providecommand{\highlight}[1]{{\color{\highlightcolor}#1}}

    \providecommand{\memotempa}{}
    \providecommand{\memotempb}{}

    \providecommand{\memo}[1]{%
      \ensuremath{%
        \framebox{\tiny\textbf{\kern-2pt\textsf{#1}}\kern-2pt}%
      }%
      \xspace}

    \RequirePackage{alphalph}
    \provideboolean{shownotes}
    \setboolean{shownotes}{true}%
    \newcounter{margnote}[page]
    \providecommand{\mgcolor}{a}
    \providecommand{\mgcolorset}[1]{\renewcommand{\mgcolor}{\alphalph{#1}}}
    \providecommand{\mgcolorsetbycounter}[1]{%
      \ifthenelse{\value{#1}<11}{%
        \renewcommand{\mgcolor}{\alph{#1}}%
      }{%
        \renewcommand{\mgcolor}{a}}%
    }
    \providecommand{\mgcolormake}{\mgcolorsetbycounter{margnote}}
    \providecommand{\mgcolorstepby}[1]{
      \setcounter{tmpcounter}{\value{margnote}}%
      \addtocounter{tmpcounter}{#1}%
      \mgcolorsetbycounter{tmpcounter}%
    }%
    \providecommand{\margnotecolor}{%
      \ifthenelse{\value{margnote}=0}{%
        \mgcolorset{10}
      }{%
        \ifthenelse{\value{margnote}<7}{%
          \mgcolormake%
        }{%
          \ifthenelse{\value{margnote}=7}{\mgcolorset{10}}{%
            \ifthenelse{\value{margnote}<11}{\mgcolormake}{%
              \ifthenelse{\value{margnote}<17}{\mgcolorstepby{-10}}{%
                \mgcolorset{10}%
              }%
            }%
          }%
        }%
      }%
    }%
    \providecommand{\margnotemark}{{\colorbox{\mgcolor}{\tiny\color{g}\upshape\texttt{\arabic{page}.\arabic{margnote}}}}}
    \providecommand{\margnote}[2][]{%
      \ifthenelse{%
        \boolean{shownotes}%
      }{%
        \stepcounter{margnote}%
        \margnotecolor%
        \margnotemark %
        \marginpar{%
          \color{\mgcolor}%
          \texttt{%
            \begin{minipage}{2cm}%
              \raggedright\tiny%
              \margnotemark%
              #2%
              \\
              {\ifx|#1|{}\else{ - #1}\fi}%
            \end{minipage}%
          }%
        }%
      }{%
      }%
    }%
    \providecommand{\mathnote}[2][]{%
      \ifthenelse{%
        \boolean{shownotes}%
      }{%
        \stepcounter{margnote}%
        \margnotecolor%
        \text{%
          \colorbox{\mgcolor}{%
            \color{g}%
            \texttt{%
              \tiny%
                  \margnotemark: %
                  \ifx|#1|{}\else{#1:}\fi%
                  #2%
            }%
          }%
        }%
      }{%
      }%
    }%
    \providecommand{\textnote}[2][]{%
      \ifthenelse{%
        \boolean{shownotes}%
      }{%
        \stepcounter{margnote}%
        \margnotecolor%
        \ \\
        \text{%
          \colorbox{\mgcolor}{%
            \begin{minipage}{.9\textwidth}
            \color{g}%
            \texttt{%
              \margnotemark: %
              \ifx|#1|{}\else{#1: }\fi%
              #2%
            }%
            \end{minipage}
          }%
        }%
      }{%
      }%
    }%

    \providecommand{\Todo}[1]{
      \ifthenelse{\boolean{shownotes}}{
        \begin{center}
        \begin{tikzpicture}
         \node[fill=a!17]{
           \begin{minipage}{\textwidth}
             \texttt{To do:}
             \\
             \texttt{\bfseries{\small #1}}
           \end{minipage}
         };
        \end{tikzpicture}
        \end{center}
      }{}}
    \provideboolean{showrevisions}
    \setboolean{showrevisions}{true}
    \newcommand{\revisionsheader}{***\newline\Warning{the following part is under revision}}
    \newcommand{\revisionsfooter}{***\newline\Warning{end of part under revision}}

    \providecommand{\Warning}[1]{    
      \begin{tikzpicture}
        \node[fill=a!27]{
          \begin{minipage}{\textwidth}
            \texttt{\bfseries{\small Warning: #1}}
          \end{minipage}
        };
      \end{tikzpicture}
    }

    \provideboolean{showcomments}
    \providecommand{\margincomment}[1]{
    \ifthenelse{\boolean{showcomments}}{\marginpar{\tiny #1}}{}
    }
    \provideboolean{showchanges}
    \setboolean{showchanges}{false}
    \providecommand{\changes}[2][]{%
      \ifthenelse{\boolean{showchanges}}{{\ifx|#1|{}\else\margnote{#1}\fi\highlight{#2}}}{#2}}
    \providecommand{\mathchanges}[2][]{%
      \ifthenelse{\boolean{showchanges}}{{\ifx|#1|{}\else\mathnote{#1}\fi\highlight{#2}}}{#2}}

    \providecommand{\changefromto}[3][replace with]{%
      \ifthenelse{\boolean{showchanges}}{%
        {\crossout{#2}\margnote{#1}}{\highlight{#3}}}{%
        #3\xspace}%
    }
    \providecommand{\ChangePar}[3][]{%
      \ifthenelse{\boolean{showchanges}}{
        {\par\textcolor{i!20}{#2}\ifx|#1|\else\margnote{#1}\fi}{\par\textcolor{a}{#3}}
      }{%
        \par #3%
      }%
    }
    \providecommand{\InsertPar}[1]{
      \ifthenelse{\boolean{showchanges}}
      {{\par$\mapsto$ \textcolor{blue}{#1}}}
      {\par #1}
    }
    
    \providecommand{\mathchangefromto}[3][]{\crossout{#2}\ifx|#1|\else\mathnote{#1}\fi\highlight{#3}}

    \providecommand{\mathscript}
    	   {\mathscr}

     \providecommand{\cA}{\ensuremath{\mathscript A}\xspace}

     \providecommand{\cE}{\ensuremath{\mathscript E}\xspace}

     \providecommand{\cH}{\ensuremath{\mathscript H}\xspace}
     \providecommand{\cI}{\ensuremath{\mathscript I}\xspace}

     \providecommand{\cW}{\ensuremath{\mathscript W}\xspace}
     \providecommand{\cX}{\ensuremath{\mathscript X}\xspace}
     \providecommand{\cY}{\ensuremath{\mathscript Y}\xspace}
     \providecommand{\cZ}{\ensuremath{\mathscript Z}\xspace}
     \providecommand{\bbbold}{\mathbb}

     \providecommand{\rN}{\ensuremath{\bbbold N}\xspace}
     
     \providecommand{\rP}{\ensuremath{\bbbold P}\xspace}
     
     \providecommand{\rR}{\ensuremath{\bbbold R}\xspace}
     
     \providecommand{\rT}{\ensuremath{\bbbold T}\xspace}

     \providecommand{\rX}{\ensuremath{\bbbold X}\xspace}
     \providecommand{\rY}{\ensuremath{\bbbold Y}\xspace}

    \providecommand{\Ae}[1][]{\ensuremath{\ifx|#1|{\ }\else{\:#1\text{-}}\fi\text{almost everywhere }}\xspace}
    \providecommand{\Aa}[1][]{\ensuremath{\text{ for }\ifx|#1|{}\else{\:#1\text{-}}\fi\text{almost all }}}
    \providecommand{\as}[1][]{\ensuremath{\ifx|#1|{\ }\else{#1\text{-}}\fi\text{almost surely}}\xspace}
    \providecommand{\aposteriori}{aposteriori\xspace}
    \providecommand{\Aposteriori}{{Aposteriori}\xspace}

     \providecommand{\naturals}{\ensuremath{\rN}}
     
     \providecommand{\NO}[1][]{\ensuremath{\naturals_0\ifx|#1|{}\else^{#1}\fi}}

     \providecommand{\reals}{\rR}

     \providecommand{\R}[1]{\reals^{#1}}
     
     \providecommand{\fieldmats}[3][F]{\csname#1\endcsname{#2\times#3}}
     \providecommand{\realmats}[2]{\fieldmats[R]{#1}{#2}}
     
     \providecommand{\RO}[1][]{{\reals_{0+}\ifx|#1|{}\else^{#1}\fi}}
     \providecommand{\RP}[1][]{{\reals_+\ifx|#1|\else^{#1}\fi}}

     \providecommand{\torus}[1]{\rT\ifthenelse{\equal{#1}1}{}{^#1}}
    
     \providecommand{\one}{\ensuremath{\bbbold 1}\xspace}
     \providecommand{\charfun}[1]{\one_{#1}}

     \providecommand{\diracdelta}[1][]{\ensuremath{{\mathrm{\delta}}\ifx|#1|{}\else_{#1}\fi}}

     \providecommand{\pic}{\ensuremath{\mathrm\pi}}
     \providecommand{\pifracl}[2][]{\fracl{\ifx|#1|\else#1\fi\pic}{#2}}
     \providecommand{\pifrac}[2][]{\frac{\ifx|#1|\else#1\fi\pic}{#2}}

     \providecommand{\closure}[1]{\overline{#1}}
     \providecommand{\inner}{\cdot}
     \providecommand{\outerp}{\wedge}

     \providecommand{\W}{\ensuremath{\varOmega}\xspace}

     \providecommand{\epsi}{\ensuremath{\epsilon}\xspace}
     \renewcommand{\epsi}{\ensuremath{\epsilon}\xspace}

     \providecommand{\qp}[2][]{\ifx|#1|\left(\else\csname#1\endcsname(\fi{#2}\ifx|#1|\right)\else\csname#1\endcsname)\fi}

     \providecommand{\qpreg}[1]{\ensuremath{(#1)}}
     \providecommand{\qpbig}[1]{\qp[big]{#1}}%
     \providecommand{\qpBig}[1]{\ensuremath{\Big(#1\Big)}}
     \providecommand{\qpbigg}[1]{\ensuremath{\bigg(\!#1\!\bigg)}}
     \providecommand{\qpBigg}[1]{\ensuremath{\Bigg(\!#1\!\Bigg)}}
     \providecommand{\qb}[2][]{\ifx|#1|\left[\else\csname#1\endcsname[\fi{#2}\ifx|#1|\right]\else\csname#1\endcsname]\fi}
     \providecommand{\qbbig}[1]{\qb[big]{#1}}%
     \providecommand{\qc}[2][]{\ifx|#1|\left\{\else\csname#1\endcsname\{\fi{#2}\ifx|#1|\right\}\else\csname#1\endcsname\}\fi}

     \providecommand{\qa}[1]{\ensuremath{\left\langle{#1}\right\rangle}}
     \providecommand{\qareg}[1]{\ensuremath{\langle#1\rangle}}
     \providecommand{\qabig}[1]{\ensuremath{\big\langle#1\big\rangle}}
     \providecommand{\qaBig}[1]{\ensuremath{\Big\langle#1\Big\rangle}}
     \providecommand{\qabigg}[1]{\ensuremath{\bigg\langle#1\bigg\rangle}}
     \providecommand{\qaBigg}[1]{\ensuremath{\Bigg\langle#1\Bigg\rangle}}

     \providecommand{\opclinter}[2]{\ensuremath{\left(#1,#2\right]}\xspace}

     \providecommand{\compowqp}[2]{\ensuremath{\qp{\!#2\!\!}^{\kern -.4em #1}\!}}
     
     \providecommand{\powqpreg}[2]{\ensuremath{%
         \qpreg{#2}^{\kern 0em\lower .1ex\hbox{\scriptsize $#1$}}\kern-.3em}}
     \providecommand{\powqpbig}[2]{\ensuremath{%
         \qpbig{#2}^{\kern -.2em\lower .3ex\hbox{\scriptsize $#1$}}\kern-.3em}}
     \providecommand{\powqpBig}[2]{\ensuremath{%
         \qpBig{#2}^{\kern -.2em\lower .3ex\hbox{\scriptsize $#1$}}\kern-.3em}}
     \providecommand{\powqpbigg}[2]{\ensuremath{%
         \qpbigg{#2}^{\kern -.2em\lower .3ex\hbox{\scriptsize $#1$}}\kern-.3em}}
     \providecommand{\powqpBigg}[2]{\ensuremath{%
         \qpBigg{#2}^{\kern -.2em\lower .3ex\hbox{\scriptsize $#1$}}}}
     \providecommand{\powp}[3][]{{#3}\ifx|#1|^{#2}\else{#1}^{#2}\fi}%
     \providecommand{\pow}[2][]{\ifx|#1|\operatorname{pow}^{#2}\else\powp{#2}{#1}\fi}%

     \providecommand{\powqp}[3][]{\powp[#1]{#2}{\qp{#3}}}%
     \providecommand{\norm}[2][]{\ifx|#1|\left|\else\csname#1\endcsname|\fi#2\ifx|#1|\right|\else\csname#1\endcsname|\fi}
     \providecommand{\normon}[2]{\norm{#1}_{#2}}

     \providecommand{\abs}[2][]{\ensuremath{\ifx|#1|{\left|#2\right|}\else{\csname#1\endcsname|{#2}\csname#1\endcsname|}\fi}}

     \providecommand{\Norm}[2][]{\ifx|#1|\left\|\else\csname#1\endcsname\|\fi{#2}\ifx|#1|\right\|\else\csname#1\endcsname\|\fi}
     \providecommand{\Normreg}[1]{\ensuremath{\|#1\|}}

     \providecommand{\Normon}[2]{\Norm{#1}_{#2}}

     \providecommand{\normonsob}[4][]{\normon{#2}{\sob{#3}{#4}\if|#1|{}\else(#1)\fi}}
     \providecommand{\Normonsob}[4][]{\Normon{#2}{\sob{#3}{#4}\if|#1|{}\else(#1)\fi}}
     \providecommand{\Normonleb}[3][]{\Normon{#2}{\leb{#3}\if|#1|{}\else(#1)\fi}}
     \providecommand{\ltwop}[3][]{\ensuremath{\qa{#2,#3}\ifx|#1|\else_{#1}\fi}}
     \providecommand{\ltwopreg}[2]{\ensuremath{\qareg{#1,#2}\ifx|#1|\else_{#1}\fi}}
     \providecommand{\ltwopbig}[2]{\ensuremath{\qabig{#1,#2}\ifx|#1|\else_{#1}\fi}}
     \providecommand{\ltwopBig}[2]{\ensuremath{\qaBig{#1,#2}\ifx|#1|\else_{#1}\fi}}
     \providecommand{\ltwopbigg}[2]{\ensuremath{\qabigg{#1,#2}\ifx|#1|\else_{#1}\fi}}
     \providecommand{\ltwopBigg}[2]{\ensuremath{\qaBigg{#1,#2}\ifx|#1|\else_{#1}\fi}}

     \providecommand{\duality}[3][]{\ensuremath{#1\langle #2\,#1\vert\,#3#1\rangle}}
     
     \providecommand{\average}[2][]{{\qa{#2}\ifx|#1|\else_{#1}\fi}}

     \providecommand{\ensemble}[2]{\ensuremath{\left\{ #1:\;#2 \right\}}}
     \providecommand{\setof}[1]{{\qc{#1}}}

     \providecommand{\conditionalto}[1]{{\left|{#1}\right.}}

    \providecommand{\measure}[1]{\ensuremath{\mathcalbf{\MakeUppercase{#1}}}}
    
    \providecommand{\probmeasure}[2][]{{\measure{#2}}\ifx|#1|\else_{#1}\fi}
    \providecommand{\Prob}{}
    \renewcommand{\Prob}[1][]{\probmeasure[{#1}]{p}}

    \providecommand{\randvars}[1][\Prob]{\operatorname{RV}\ifx|#1|{}\else{(#1)}\fi}
    \providecommand{\discrandvars}[1][\Prob]{\operatorname{DRV}\ifx|#1|{}\else{({#1)}\fi}} 
    \providecommand{\contrandvars}[1][\Prob]{\ensuremath{\operatorname{CDRV}\ifx|#1|{}\else(#1)\fi}} 
     \def\env@matrix{\hskip -\arraycolsep
      \let\@ifnextchar\new@ifnextchar
      \array{*\c@MaxMatrixCols c}}
     \makeatletter
     \renewcommand*\env@matrix[1][c]{\hskip -\arraycolsep
       \let\@ifnextchar\new@ifnextchar
       \array{*\c@MaxMatrixCols #1}}
     \makeatother
     \providecommand{\irow}[2]{#1_{#2}}%
     \providecommand{\icol}[2]{#1^{#2}}%

     \providecommand{\ijrowcol}[3]{\icol{\irow{#1}{#2}}{#3}}
     \providecommand{\entry}[1]{\qb{#1}}
     \providecommand{\vecentry}[2]{\irow{#1}{#2}}

     \providecommand{\rowof}[1]{\qb{#1}}

     \providecommand{\getentryi}[2]{\irow{\entry{#1}}{#2}}
     
     \providecommand{\getvecentry}[2]{\getentryi{\vec #1}{#2}}

     \providecommand{\dismatof}[2][r]{\begin{bmatrix}[#1]#2\end{bmatrix}}

     \providecommand{\matentry}[3]{\ijrowcol{#1}{#2}{#3}}

     \providecommand{\block}[5]{\ijrowcol{#1}{\ifx#2#3{\rowof{#2}}\else\rowof{{#2}\dotsc{#3}}\fi}{\ifx#4#5{\rowof{#4}}\else\rowof{{#4}\dotsc{#5}}\fi}}
     \providecommand{\colblock}[3]{\getvecentry{#1}{\ifx#2#3{#2}\else\fromto{#2}{#3}\fi}}

     \providecommand{\dismatskeldots}[4]{
       \dismatof[c]{
         #1&\dotsc&#3
         \\
         \vdots & \ddots &\vdots
         \\
         #2&\dotsc&#4
       }
     }
     \providecommand{\dismatcommfromtofromto}[5]{
       \dismatskeldots{#1#2#4}{#1#3#4}{#1#2#5}{#1#3#5}
     }
     \providecommand{\dismatcustfromtofromto}[6][matentry]{
       \dismatcommfromtofromto{\csname#1\endcsname{#2}}#3#4#5#6
     }
     \providecommand{\dismatcustfromtofromto}[6][matentry]{
       \dismatskeldots{%
         \csname#1\endcsname{#2}{#3}{#4}%
       }{%
         \csname#1\endcsname{#2}{#3}{#6}%
       }{%
         \csname#1\endcsname{#2}{#5}{#4}%
       }{%
         \csname#1\endcsname{#2}{#5}{#6}%
       }%
     }%
     \providecommand{\dismatcustfromtofromto}[6][matentry]{
       \dismatof{
         \csname#1\endcsname{#2}{#3}{#4}&\dotsc&\csname#1\endcsname{#2}{#3}{#6}
         \\
         \vdots & \ddots &\vdots
         \\
         \csname#1\endcsname{#2}{#5}{#4}&\dotsc&\csname#1\endcsname{#2}{#5}{#6}
       }
     }

     \providecommand{\dissysaxbdotsnm}[5]{\begin{matrix}[r]%
         \matentry{#1}11\vecentry{#2}1&+\dotsb&+\matentry{#1}1{#5}\vecentry{#2}{#5}
         &
         =
         \ifx|#3|0\else{\vecentry {#3}1}\fi
         \\
         \dotsb
         \\
         \matentry{#1}{#4}1\vecentry{#2}1&+\dotsb&+\matentry{#1}{#4}{#5}\vecentry{#2}{#5}
         &
         =
         \ifx|#3|0\else{\vecentry {#3}{#4}}\fi
     \end{matrix}}
     \providecommand{\seqof}[1]{\qp{#1}}%
     \providecommand{\seqs}[2]{\seqof{#1}_{#2}}
     \providecommand{\sets}[2]{\setof{#1}_{#2}}%
     \providecommand{\seqi}[3][]{\seqs{#2_{#3}}{\ifx|#1|{#3}\else{{#3}\in{#1}}\fi}}%

     \providecommand{\seti}[3][]{\sets{#2_{#3}}{\ifx|#1|_{#3}\else_{{#3}\in{#1}}\fi}}%
     \providecommand{\sequ}[3][]{\seqs{#2^{#3}}{\ifx|#1|{#3}\else{{#3}\in{#1}}\fi}}%
     \providecommand{\setu}[3][]{\sets{#2^{#3}}{\ifx|#1|{#3}\else{{#3}\in{#1}}\fi}}%
     \providecommand{\limofat}[3][]{\ensuremath{\lim_{\ifx|#1|{}\else{#1\ni}\fi#3}{#2}}}
     \providecommand{\limsupofat}[3][]{\ensuremath{\limsup_{\ifx|#1|{}\else{#1\ni}\fi#3}{#2}}}
     \providecommand{\liminfofat}[3][]{\ensuremath{\liminf_{\ifx|#1|{}\else{#1\ni}\fi#3}{#2}}}

     \providecommand{\listdotsfrom}[3][]{\ensuremath{#2\ifx|#1|\else#1\fi,#3\ifx|#1|\else#1\fi,\dotsc}}
     \providecommand{\listdotsfromto}[3][]{\ensuremath{#2\ifx|#1|\else#1\fi,\dotsc,#3\ifx|#1|\else#1\fi}}
     \providecommand{\listifromto}[5][]{\ensuremath{{#2}_{#3}\ifx|#1|\else#1\fi},\text{ for }\ensuremath{\rangefromto{#3}{#4}{#5}}\xspace}
     \providecommand{\listufromto}[5][]{\ensuremath{{#2}^{#3}\ifx|#1|\else#1\fi},\text{ for }\ensuremath{\rangefromto{#3}{#4}{#5}}\xspace}
     \providecommand{\listitwo}[2][]{\ensuremath{#2_1\ifx|#1|\else#1\fi,#2_2\ifx|#1|\else#1\fi}}
     \providecommand{\listutwo}[2][]{\ensuremath{#2^1\ifx|#1|\else#1\fi,#2^2\ifx|#1|\else#1\fi}}
     \providecommand{\listithree}[2][]{\ensuremath{#2_1\ifx|#1|\else#1\fi,#2_2\ifx|#1|\else#1\fi,#2_3\ifx|#1|\else#1\fi}}
     \providecommand{\listuthree}[2][]{\ensuremath{#2^1\ifx|#1|\else#1\fi,#2^2\ifx|#1|\else#1\fi,#2^3\ifx|#1|\else#1\fi}}

     \providecommand{\listidotsfromto}[4][]{\listdotsfromto[#1]{#2_{#3}}{#2_{#4}}}

     \providecommand{\sumifromto}[3]{\ensuremath{\sum_{#1=#2}^{#3}}}

     \providecommand{\jump}[2][]{\ensuremath{\left\llbracket #2\right\rrbracket\ifx|#1|{}\else_{#1}\fi}}

     \providecommand{\fromto}[2]{\ensuremath{\setof{#1\dotsc#2}}}%

     \providecommand{\integerbetween}[2]{\ensuremath{={#1},\dotsc,{#2}}}
     
     \providecommand{\rangefromto}[3]{\ensuremath{#1\integerbetween{#2}{#3}}}

     \providecommand{\d}{}
     \renewcommand{\d}[1][]{\ensuremath{\operatorname{d}\!\ifx|#1|\else{_{#1}}\fi}}
     
     \providecommand{\ds}[1][]{\d{\measure S}}
     \providecommand{\D}[1][]{\ensuremath{\operatorname{D}\!\ifx|#1|\else{_{#1}}\fi}}

    \providecommand{\registered}%
    {\ensuremath{^\text{\textregistered}}}

    \providecommand{\AND}{\ensuremath{\text{ and }}}%
    \providecommand{\tand}{\ensuremath{\text{ and }}}

    \providecommand{\constant}[1]{\ensuremath{C_{#1}}}
    \providecommand{\constext}[2][]{\constant{\textup{#2}{\ifx|#1|{}\else{,\ensuremath{#1}}\fi}}}            %
    \providecommand{\constref}[2][]{\ensuremath{\constant{\textup{\ref{#2}{\ifx|#1|{}\else{,\ensuremath{#1}}\fi}}}}}
    \providecommand{\constdef}[2][]{\label{#2}\ensuremath{\constant{\textup{\ref{#2}{\ifx|#1|{}\else{,\ensuremath{#1}}\fi}}}}}

    \providecommand{\funkref}[3][]{\ensuremath{{#3}_{\textup{\ref{#2}{\ifx|#1|{}\else{,\ensuremath{#1}}\fi}}}}}

    \renewcommand{\div}[1][]{\nabla\ifx|#1|{}\else\kern-2pt_{#1}\fi\kern-2pt\inner}
    \providecommand{\divof}[2][]{\div[#1]\ifx|#2|{}\else\qb{#2}\fi}
    \providecommand{\divideabyb}[2]{\operatorname{div}(a,b)}
    \providecommand{\grad}{}
    \renewcommand{\grad}[1][]{\nabla\ifx|#1|\else_{#1}\fi}

    \providecommand{\rot}[1][]{\nabla\ifx|#1|\else_{#1}\fi\outerp}
    
    \providecommand{\rowdiv}[1][]{\D\ifx|#1|{}\else\kern-1pt_{#1}\kern-2pt\fi\cdot}
    \providecommand{\rowdivof}[2][]{\rowdiv[#1]\ifx|#2|{}\else\qb{#2}\fi}

    \providecommand{\esssup}{\operatorname{ess\,sup}}         %

    \providecommand{\inv}[1][]{\operatorname{inv}\ifx|#1|\else^{#1}\fi}
    \providecommand{\ivt}[1]{\operatorname{ivt}\ifx|#1|\else^{#1}\fi}
    \providecommand\tensorinvariant\ivt

    \providecommand{\mod}{}
    \renewcommand{\mod}[1][]{\operatorname{mod}\ifx|#1|\else\kern-1pt_{#1}\fi}
    
    \providecommand{\fracl}[3][]{\ifx|#1|\nicefrac{#2}{#3}\else{#2}#1/{#3}\fi}

    \providecommand{\qpfracl}[3][]{\qp{\ifx|#1|\fracl{#2}{#3}\else{#2}#1/{#3}\fi}}
    \providecommand{\qpfrac}[3][]{\qp{\ifx|#1|\frac{#2}{#3}\else{#2}#1/{#3}\fi}}
    \providecommand{\absfracl}[3][]{\abs{\ifx|#1|\fracl{#2}{#3}\else{#2}#1/{#3}\fi}}
    \providecommand{\absfrac}[3][]{\abs{\ifx|#1|\frac{#2}{#3}\else{#2}#1/{#3}\fi}}
    \providecommand{\fraclff}[3][]{\ifx|#1|{#2}/{#3}\else{#2}#1/{#3}\fi}

    \providecommand{\eye}[1][]{\vec{\mathrm I}\ifx|#1|{}\else_{#1}\fi}%
    \providecommand{\numeye}[1][]{\boldsymbol{\mathsf{I}}\ifx|#1|{}\else_{#1}\fi}%
    \providecommand{\Eye}[1]{
      \begin{bmatrix}
      \ifthenelse{#1>1}{
        \ifthenelse{#1>2}{
          \ifthenelse{#1>3}{
            \ifthenelse{#1>4}{
              1&\zeroentry&\dotso&\zeroentry
              \\
              \zeroentry&1&\dotso&\zeroentry
              \\
              \vdots&\vdots&\ddots&\vdots
              \\
              \zeroentry&\zeroentry&\dotso&1
            }{        
              1&\zeroentry&\zeroentry&\zeroentry
              \\
              \zeroentry&1&\zeroentry&\zeroentry
              \\
              \zeroentry&\zeroentry&1&\zeroentry
              \\
              \zeroentry&\zeroentry&\zeroentry&1
            }
          }{
            1&\zeroentry&\zeroentry
            \\
            \zeroentry&1&\zeroentry
            \\
            \zeroentry&\zeroentry&1
          }
        }{
          1&\zeroentry
          \\
          \zeroentry&1
        }
      }{
        1
      }
      \end{bmatrix}
    }

    \providecommand{\lebmeas}[1][]{\measure L^{#1}}     %
    \providecommand{\lebmeasof}[2][]{\ifx|#1|\left|#2\right|\else\lebmeas[#1]\qp{#2}\fi}         %
    \providecommand{\meshsize}[1][]{h\ifx|#1|\else_{#1}\fi}
    \providecommand{\maxi}[2]{#1\vee#2}                       %

    \providecommand{\mini}[2]{#1\wedge#2}                     %
    \providecommand{\dash}[1][']{\ifthenelse{\equal{#1}{'}\OR\equal{#1}{''}}{#1}{^{(#1)}}}
    
    \providecommand{\pdfrac}[2][]{\ensuremath{\frac{\partial\ifx|#1|\phantom{#2}\else{#1}\fi}{\partial{#2}}}} %
    \providecommand{\pdfracpow}[3][]{\ensuremath{\frac{\partial^{#3}\ifx|#1|\phantom{#2}\else{#1}\fi}{\partial{#2}^{#3}}}} %

    \providecommand{\pd}[2][]{\ensuremath{\partial_{#2}}{\ifx|#1|{}\else{\qb{#1}}\fi}} %

    \providecommand{\pdt}[1][]{\pd[#1]t}                       %

    \renewcommand{\Im}{\operatorname{im}}                 %
    \renewcommand{\Re}{\operatorname{re}}                 %
    \providecommand{\imaginpart}[1][]{\Im{\ifx|#1|{}\else\qp{#1}\fi}} %
    \providecommand{\realpart}[1][]{\Re{\ifx|#1|{}\else\qp{#1}\fi}} %
    \providecommand\determinant\det

    \providecommand{\transpose}{\intercal}%

    \providecommand{\Transpose}[1]{\ensuremath{{#1}^{\transpose}}}
    
    \providecommand{\Transposevec}[1]{\Transpose{\vec{#1}}}
    \providecommand{\transposevec}[1]{\Transposevec{#1}}

    \providecommand{\orthogonalto}[1][]{\ensuremath{\perp\ifx|#1|{}\else{_{#1}}\fi}}
    
    \providecommand{\rowof}[1]{\ensuremath{\vecof{#1}}}

    \providecommand{\zeroentry}{\phantom0}%

    \providecommand{\smint}{\ensuremath{{\text{\textbf{/}}}\kern-.75em\smallint}}
    \renewcommand{\smint}[1][]{\lower12.3pt\hbox{\begin{tikzpicture}\draw[line width=.75pt] (-3pt,-0.5)--(1pt,-0.5) node[pos=0.6]{$\int$};\path (3pt,-24pt)node {\scriptsize $#1$};\end{tikzpicture}}}

    \providecommand{\lap}{\ensuremath{\Delta}}
    \providecommand{\lapin}[1][]{\lap\ifx|#1|\else_{#1}\fi}
    \providecommand{\normalsymbol}{\operatorname{n}}
    \providecommand{\normal}[1][]{\normalsymbol\ifx|#1|\else_{#1}\fi}%
    \providecommand{\normalto}[2][]{\ensuremath{\normal[#2]\ifx|#1|\else\qp{#1}\fi}}
    \providecommand{\normalder}[1][]{\ensuremath{\normal\ifx|#1|\else\qp{#1}\fi{\inner\grad}}}

    \providecommand{\tangentialsymbol}{\operatorname{t}}
    \providecommand{\tangentialto}[2][]{\tangentialsymbol\ifx|#1|\else^{#1}\fi\ifx|#2|\else_{#2}\fi}

    \providecommand{\intersected}{\ensuremath{\cap}}
    
    \providecommand{\meet}{\intersected}

    \providecommand{\union}[1]{\ensuremath{\bigcup\nolimits_{#1}}}

    \providecommand{\unions}[3][]{\union{#2\in{#3}\ifx|#1|\else:#1\fi}}

    \ifthenelse{\boolean{usesvjour}}{}{
      \let\vec\undefined
      \providecommand{\vec}[1]{\ensuremath{\boldsymbol{#1}}}
      \renewcommand{\vec}[1]{\ensuremath{\boldsymbol{#1}}}
    }

    \providecommand{\hatmat}[1]{\hat{\mat{#1}}}

    \providecommand{\geomat}[1]{\vec{\MakeUppercase{#1}}}

    \providecommand{\mat}[1]{\geomat{#1}} %
    \providecommand{\Prob}[1][]{\ensuremath{\operatorname{Prob}\ifx|#1|{}\else_{#1}\fi}}

    \providecommand{\pdf}[2][]{\ensuremath{\operatorname{pdf}_{#2\ifx|#1|{}\else{\conditionalto{#1}}\fi}}\xspace}

    \providecommand{\expectation}{\ensuremath{\operatorname{E}}}
    \providecommand{\EX}[1][]{\ensuremath{\expectation\ifx|#1|{}\else_{#1}\fi}}

    \providecommand{\gausskernel}[3][x]{%
      \ensuremath{
        \exp\frac{-\if#20{#1}\else(#1-\mu)\fi^2}{%
          2\if#31{}\else\powp2{#3}\fi}%
      }%
    }
    \providecommand{\gaussdistribution}[3][x]{%
      \ensuremath{\frac1{\sqrt{2\pic}\if#31{}\else#3\fi}%
        \gausskernel[#1]{#2}{#3}
      }%
    }%

    \providecommand{\boundary}{\partial}

    \providecommand{\SPD}{\operatorname{SPD}}
    \providecommand{\spdmats}[2][F]{\SPD(\csname#1\endcsname{#2})}
     \providecommand{\Continuous}{\ensuremath{\operatorname C}\xspace}%
     \providecommand{\Hspace}{\ensuremath{\operatorname H}\xspace}
     \providecommand{\Lebesgue}{\ensuremath{\operatorname L}\xspace}
     \providecommand{\Besovspace}{\ensuremath{\operatorname B}\xspace}

     \providecommand{\Weaklyder}{\ensuremath{\operatorname W}\xspace}
     
     \providecommand{\dual}[1]{\ensuremath{{#1}'}}
     
     \providecommand{\dualspace}[2][]{\dual{\linspace{#2}\ifx|#1|\else{_{#1}}\fi}}
     \providecommand{\bidual}[1]{\ensuremath{{#1}''}}
     \providecommand{\bidualspace}[2][]{\bidual{\linspace{#2}\ifx|#1|\else{_{#1}}\fi}}

     \providecommand{\cont}[1]{\ensuremath{\Continuous^{#1}}}

     \providecommand{\BV}[1]{\ensuremath{\operatorname{BV}}}

     \providecommand{\leb}[1]{\ensuremath{\Lebesgue_{#1}}}
     \providecommand{\lebloc}[1]{\ensuremath{{{\Lebesgue}^{\kern-.20em\lower .1ex\hbox{\tiny\textrm{\textup{loc}}}}_{#1}}}}
     \providecommand{\lebnorm}[3][]{\ensuremath{\Norm{#2}_{\leb{#3}\ifx|#1|{}\else(#1)\fi}}}
     \providecommand{\bes}[3][]{\ensuremath{\Besovspace^{#2}_{#3\ifx|#1|\else,#1\fi}}}

     \providecommand{\sob}[2]{\ensuremath{{\smash\Weaklyder}^{#1}_{#2}}}

     \providecommand{\sobh}[1]{\ensuremath{\Hspace^{#1}}}
     \providecommand{\vecsobh}[1]{\ensuremath{\vec\Hspace^{#1}}}
     \providecommand{\hdiv}[1][]{\vecsobh{\operatorname{div}}\ifx|#1|\else(#1)\fi}
     \providecommand{\hcurl}[1][]{\vecsobh{\operatorname{curl}}\ifx|#1|\else(#1)\fi}
     
     \providecommand{\sobhz}[2][]{\sobh{#2}_{0\ifx+#1+\else|#1\fi}}

     \providecommand{\Lip}[1][]{\ensuremath{\operatorname{Lip}}\ifx|#1|{}\else{\qp{#1}}\fi}

     \providecommand{\poly}[1]{\ensuremath{\rP}^{#1}}

     \providecommand{\Symmatrices}[2][R]{\ensuremath{\operatorname{Sym}{(\csname#1\endcsname{#2})}}}
     
     \providecommand{\SAmatrices}[2][F]{\ensuremath{\operatorname{SA}{(\csname#1\endcsname{#2})}}}
     \providecommand{\mesh}[2][]{{\ensuremath{\mathcalbf{\MakeUppercase{#2}}\ifx|#1|\else_{#1}\fi}}}

    \providecommand{\crouzeixraviart}[1][1]{\operatorname{CR}\ifx|#1|{}\else{^{#1}}\fi}

    \providecommand{\linspace}[1]{\mathscript{\MakeUppercase{#1}}}

    \providecommand{\clinopss}[2]{\clinopss{\linspace{#1}}{\linspace{#2}}}
    \providecommand{\fepartition}[2][]{\mathscript{\MakeUppercase{#2}}\ifx|#1|{}\else_{#1}\fi}
    \providecommand{\fespace}[2][]{\mathbb{\uppercase{#2}}\ifx|#1|{}\else_{#1}\fi}
    
    \providecommand{\vespace}[1][]{\fespace v\ifx|#1|\else_{#1}\fi}

    \providecommand{\fes}[2]{\ensuremath{\fespace{#1}^{#2}}}

    \providecommand{\fe}[2][]{\ensuremath{\uppercase{#2}\ifx|#1|\else_{#1}\fi}}%

    \providecommand{\vecfe}[2][]{\ensuremath{\vec{\uppercase{#2}}\ifx|#1|{}\else{_{#1}}\fi}}%
    
    \providecommand{\matfe}[2][]{\ensuremath{\mat{\uppercase{#2}}\ifx|#1|{}\else{_{#1}}\fi}}%
    
    \providecommand{\hatmatfe}[2][]{\ensuremath{\hatmat{\uppercase{#2}}\ifx|#1|{}\else{_{#1}}\fi}}%

    \providecommand{\Foreach}{\text{ for each }}%

    \RequirePackage{amscd}

    \providecommand{\funk}[3]{\ensuremath{#1:#2\to#3}}

    \providecommand{\dfunkmapsto}[6][]{\ensuremath{
        \begin{array}{rrcl}
          {#2}: & {#4} &  \to   & {#6}
          \\
                & {#3} &\mapsto & {#5\text{\ #1}}
        \end{array}\quad}}

    \providecommand{\implies}{\ensuremath{\:\Rightarrow\:}\xspace}
    \renewcommand{\implies}{\ensuremath{\:\Rightarrow\:}\xspace}

    \providecommand{\imbedded}{{\ensuremath{\,\hookrightarrow\,}}}
    
    \providecommand{\embedsin}{\imbedded}

    \providecommand{\restriction}[2]{\left.#1\right|_{#2}}
    \renewcommand{\restriction}[2]{\left.#1\right|_{#2}}

    \providecommand{\evalat}[3][]{\qb{#2}_{\ifx|#1|{}\else#1=\fi#3}}
    \providecommand{\evaldiff}[4][]{\qb{#2}^{\ifx|#1|{}\else#1=\fi#3}_{\ifx|#1|{}\else#1=\fi#4}}

    \providecommand\bs{\char '134}   %

    \providecommand{\texcommand}[1]{\texttt{\bs{\nolinkurl{#1}}}\xspace}
    
    \providecommand{\codename}[1]{\nolinkurl{#1}\xspace}

    \ifthenelse{\boolean{useutopia}}{%
      
    }{%
      
    }
    \providecommand{\indexen}[2][]{{\ifthenelse{\boolean{shownotes}}{\color b}{}#2\ifx|#1|\index{#2}\else\index{#1}\fi}}
    \providecommand{\ListParameters}{}
    \renewcommand{\ListParameters}%
    {
    	 \setlength{\topsep}{0pt}
    	 \setlength{\leftmargin}{0pt}
             \setlength{\itemsep}{0pt}
    	 \setlength{\parsep}{0pt}
    	 \setlength{\parskip}{0pt}
             \setlength{\labelsep}{0pt}
    	 \setlength{\itemindent}{0pt}
    }
    {%
      \begin{list}%
        {}%
        {\ListParameters%
        
    }}%
    {\end{list}}
    \newcounter{tmpcounter}
    \newcounter{LetterListItem}
    \renewcommand{\theLetterListItem}{(\alph{LetterListItem})}

    \newcounter{CapitalListItem}
    \renewcommand{\theCapitalListItem}{\Alph{CapitalListItem}.}

    \newcounter{NumberListItem}
    \renewcommand{\theNumberListItem}{\arabic{NumberListItem}}
    {
    	\begin{list}%
    	{\theNumberListItem.\ }%
    	{\usecounter{NumberListItem}%
    	 \ListParameters
    	}
    }%
    {\end{list}}
    \newcounter{QuestionListItem}
    \renewcommand{\theQuestionListItem}{\textbf{Question \arabic{QuestionListItem}}}
    {
    	\begin{list}%
    	{\theQuestionListItem.\ }%
    	{\usecounter{QuestionListItem}%
    	 \ListParameters
    	}
    }%
    {\end{list}}
    \newcounter{RomanListItem}
    \renewcommand{\theRomanListItem}{(\roman{RomanListItem})}
    {
    	\begin{list}%
    	{\theRomanListItem\ }%
    	{\usecounter{RomanListItem}
    	 \ListParameters
    	}
    }%
    {\end{list}}
    \newcounter{StepsItem}
    {
    	\begin{list}%
    	{Step \theStepsItem.\ }%
    	{\usecounter{StepsItem}%
    	 \ListParameters
    	}
    }%
    {\end{list}}
    \newcounter{CasesListItem}
    \renewcommand{\theCasesListItem}{\Alph{CasesListItem}}
    {
    	\begin{list}%
    	{\emph{Case \theCasesListItem.}\ }%
    	{\usecounter{CasesListItem}%
    	 \ListParameters
    	}
    }%
    {\end{list}}
    \newcounter{QAListItem}
    \renewcommand{\theQAListItem}{Q\arabic{QAListItem}:}
    {
    	\begin{list}%
    	{\theQAListItem}%
    	{\usecounter{QAListItem}
    	 \ListParameters
    	}
    }%
    {\end{list}}

    \ifthenelse{\boolean{isthesis}}{%
      \setcounter{secnumdepth}{1}%
    }{
      \setcounter{secnumdepth}{2} %
    }
    \providecommand{\ListParameters}{}
    \renewcommand{\ListParameters}
    {
    	 \setlength{\topsep}{0em}
    	 \setlength{\leftmargin}{0em}
             \setlength{\itemsep}{0ex}
    	 \setlength{\parsep}{.5ex}
    	 \setlength{\itemindent}{\labelsep}
    	 \addtolength{\itemindent}{\labelwidth}
    }

      \providecommand{\ObsName}{Remark}%
      \providecommand{\RemName}{Remark}%
      \providecommand{\NotName}{Notation}%
      \providecommand{\BFNName}{Big~Fat~Note}%
      \providecommand{\DefName}{Definition}%
      \providecommand{\ExaName}{Example}%
      \providecommand{\TheName}{Theorem}%
      \providecommand{\LemName}{Lemma}%
      \providecommand{\ProName}{Proposition}%
      \providecommand{\CorName}{Corollary}%
      \providecommand{\PbmName}{Problem}%
      \providecommand{\HypName}{Hypothesis}%
      \providecommand{\AlgName}{Algorithm}%
      \providecommand{\ExeName}{Exercise}%
      \providecommand{\SolName}{Solution}%
      \providecommand{\ClaName}{Claim}%
      \providecommand{\EsyName}{Essay}%
      \providecommand{\Proofname}{Proof}%
      \providecommand{\Derivename}{Derivation}%
    
    \ifthenelse{\boolean{isthesis}}{%
      \providecommand{\Thecounter}{The}
    }{%
      \providecommand{\Thecounter}{subsection}
    }
    \makeatletter
    \newcommand{\oltikzgetxy}[3]{%
      \tikz@scan@one@point\pgfutil@firstofone#1\relax
      \edef#2{\the\pgf@x}%
      \edef#3{\the\pgf@y}%
    }
    \makeatother
    \providecommand{\pdfformat}[1]{
       \provideboolean{pdfoutput}
       \setboolean{pdfoutput}{#1}%
      \ifthenelse{\boolean{pdfoutput}}{
        \typeout{using pdf}
\usepackage{pdfsync}
        \providecommand{\graphext}{pdf}
        \renewcommand{\graphext}{pdf}
        \providecommand{\graphextex}{pdf_t}
        \renewcommand{\graphextex}{pdf_t}
      }{
        \typeout{using eps}
        \RequirePackage[dvips]{graphicx,xcolor}
        \providecommand{\graphext}{eps}
        \renewcommand{\graphext}{eps}
        \providecommand{\graphextex}{eps_t}
        \renewcommand{\graphextex}{eps_t}
      }
      \RequirePackage{epsfig}
      \RequirePackage{tikz}
      \RequirePackage{rotating}
      \RequirePackage{graphicx}
      \RequirePackage{xcolor}
      \provideboolean{darkcolortheme}
      \definecolor{SussexFlint}{rgb}{.00,.19,.21}
      \definecolor{SussexGrey}{rgb}{.51,.58,.49}
      \definecolor{SussexOrange}{rgb}{.94,.29,.00}
      \definecolor{SussexYellow}{rgb}{1.00,.73,.00}
      \definecolor{SussexRed}{rgb}{.94,.01,.49}
      \definecolor{SussexPurple}{rgb}{.48,.06,.44}
      \definecolor{SussexGreen}{rgb}{.00,.58,.46}
      \definecolor{OmarGreen}{rgb}{.00,.68,.36}
      \definecolor{SussexBlue}{rgb}{.00,.58,.65}
      \definecolor{OmarBlue}{rgb}{.00,.38,.65}
      \colorlet{a}{OmarBlue}%
      \colorlet{b}{SussexOrange}
      \colorlet{c}{SussexGreen}
      \colorlet{d}{SussexPurple}%
      \colorlet{e}{SussexRed}
      \colorlet{f}{SussexYellow}
      \colorlet{g}{white}%
      \colorlet{h}{SussexGrey}%
      \colorlet{i}{black}%
      \colorlet{j}{SussexFlint}
      \colorlet{colora}{a}
      \colorlet{colorb}{b}
      \colorlet{colorc}{c}
      \colorlet{colord}{d}
      \colorlet{colore}{e}
      \colorlet{colorf}{f}
      \colorlet{colorg}{g}
      \colorlet{colorh}{h}
      \colorlet{colori}{i}
      \colorlet{colorj}{j}
      \newcommand{\mausDarkColorTheme}{
        \colorlet{a}{SussexYellow!50!yellow}
        \colorlet{b}{SussexBlue}%
        \colorlet{c}{SussexRed!50!red}
        \colorlet{d}{SussexOrange!50!yellow}
        \colorlet{e}{SussexGreen!50!green}
        \colorlet{f}{SussexPurple!50!magenta}
        \colorlet{g}{black}%
        \colorlet{h}{SussexFlint!50!black}
        \colorlet{i}{white}%
        \colorlet{j}{SussexGrey}
      }
      \ifthenelse{\boolean{darkcolortheme}}{\mausDarkColorTheme}{}
    }
    \providecommand{\solution}{\textbf{\SolName.}\xspace}

     \newcounter{phantombox}[enumi]%
     \provideboolean{showphantoms}
     \renewcommand{\thephantombox}{\Alph{phantombox}}%
     \newcommand{\phantombox}[1]{\stepcounter{phantombox}%
       \ensuremath{\boxed{%
           {\ifthenelse{\boolean{showphantoms}}{#1}{\phantom{#1}}}%
           {\texttt{\tiny\ \colorbox{i!50}{\color g\thephantombox}}
           }%
         }%
       }%
     }

     \provideboolean{hidesolution}
     \newcommand{\consolution}[2][]{
       \ifthenelse{\boolean{hidesolution}}{#1\setboolean{showphantoms}{false}}{%
         {\setboolean{showphantoms}{true}\color{i!50}\par \small {\solution}\ #2\par\ \\[5pt]}}
     }
     \provideboolean{showmarks}
     \providecommand{\showmarks}[1]{%
       \ifthenelse{%
         \boolean{showmarks}}{%
         \marginpar{%
           \tiny [$#1$ mark\ifthenelse{\equal{#1}1}{\phantom{s}}s]}%
       }{}}%

     \newcommand{\condibreak}{\ifthenelse{\boolean{hidesolution}}{\newpage}{}}

     \providecommand{\qeyword}[1]{\index{#1}\ifthenelse{\boolean{shownotes}}{\texttt{\color{e}[#1]}}{}}
     \providecommand{\sourcecite}[2][]{\index{#1}\ifthenelse{\boolean{shownotes}}{\texttt{\color{d}[source: \cite[#1]{#2}]}}{}}

     \RequirePackage{hyphenat}
    \RequirePackage{lineno}
    \ifthenelse{\boolean{showchanges}}{
      \newcommand{\llabel}[1]{\hypertarget{llineno:#1}{\linelabel{#1}}}
      \newcommand{\lref}[1]{\hyperlink{llineno:#1}{\ref*{#1}}}
    }{
      \newcommand\llabel[1]{}
      \newcommand\lref[1]{}
    }
    \provideboolean{includeresponses}
    \setboolean{includeresponses}{false}
    \providecommand{\mailto}[1]{\href{mailto:#1}{\nolinkurl{#1}}}

\newtheoremstyle{plain}%
  {}%
  {}%
  {\mdseries\slshape}%
  {\parindent}%
  {\bfseries}%
  {.}%
  {.5em}%
  {}%

\newtheoremstyle{note}%
  {}%
  {}%
  {}%
  {\parindent}%
  {\bfseries}%
  {.}%
  {.5em}%
  {}%

\newtheoremstyle{claim}%
  {}%
  {}%
  {\mdseries\slshape}%
  {}%
  {\bfseries}%
  {}%
  {.5em}%
  {}%

\newtheoremstyle{exercise}%
  {}%
  {}%
  {}%
  {}%
  {\bfseries}%
  {.}%
  {1em}%
  {}%

\newtheoremstyle{break}%
  {}%
  {}%
  {}%
  {}%
  {\bfseries}%
  {.}%
  {\newline}%
  {}%

\swapnumbers{
  \theoremstyle{plain}
  \ifthenelse
      {\boolean{isthesis}}
      {\newtheorem{The}{\TheName}[section]}%
      {
      }%
{
   \theoremstyle{plain}

   \renewcommand{\Thecounter}{subsection}

   \newtheorem*{The*}{\TheName}
   \newtheorem*{Lem*}{\LemName}
   \newtheorem*{Pro*}{\ProName}
   \newtheorem*{Cor*}{\CorName}
   \newtheorem*{Pbm*}{\PbmName}
   \newtheorem*{Hyp*}{\HypName}
   \newtheorem*{Exe*}{\ExeName}
   \newtheorem*{Txx*}{\ExeName} %
   \newtheorem*{Con*}{Conclusion}
   \newtheorem*{Sum*}{Summary}
 }
 {
   \theoremstyle{claim}

 }
 {
   \theoremstyle{note}

   \newtheorem*{Obs*}{\ObsName}

   \newtheorem*{Def*}{\DefName}
   \newtheorem*{Exa*}{\ExaName}
   \newtheorem*{Alg*}{\AlgName}
 }

 {
   \theoremstyle{break}
 }
}

\newenvironment{The}[1][]{%
  \ifx&#1&%
  \subsection{\TheName\xspace}%
  \else%
  \subsection[#1 theorem]{\TheName\ (#1)}%
  \fi%
  \slshape}{%
  \upshape}
\newenvironment{Pro}[1][]{\subsection{\ProName\xspace{\ifx&#1&{}\else{ (#1)}\fi}}\slshape}{\upshape}
\newenvironment{Lem}[1][]{\subsection{\LemName\xspace{\ifx&#1&{}\else{ (#1)}\fi}}\slshape}{\upshape}

\newenvironment{Obs}[1][]{\subsection{\ObsName\xspace{\ifx&#1&{}\else{ (#1)}\fi}}}{}
\newenvironment{Exa}[1][]{\subsection{\ExaName\xspace{\ifx&#1&{}\else{ (#1)}\fi}}}{}

\providecommand{\qed}{\vrule height 5pt depth 0pt width 3pt}
\providecommand{\qqed}{{\raggedright{\ \hfill\qed}}}

\newcounter{passo}

{\par\noindent{\bf \Proofname\ #1}\setcounter{passo}{0}}%
{\qqed\par}
{\par\noindent{\bf \Derivename\ #1}\setcounter{passo}{0}}%
{\qqed\par}
\newenvironment{Proof*}[1][{}]%
{\subsection{\Proofname\ #1}\setcounter{passo}{0}}
{\qqed\par}

\setboolean{usemathrsfs}{false}%
\usepackage{enumerate}
\numberwithin{equation}{section}
\newcommand{\honezw}{{\sobhz1(\W)}} %

\newcommand{\tm}{{t_m}}

\newcommand{\timestep}{\tau}

\newcommand{\ellop}{\cA}
\newcommand{\abil}[2]{\ensuremath{a\left(#1,#2\right)}}

\providecommand{\aposteriori}{{a~posteriori}\xspace}
\providecommand{\Aposteriori}{{A~posteriori}\xspace}

\renewcommand{\aposteriori}{{a~posteriori}\xspace}
\renewcommand{\Aposteriori}{{A~posteriori}\xspace}

\providecommand{\ellop}{\ensuremath{\cA}}
\providecommand{\pivot}{\ensuremath{\cH}}
\providecommand{\elldom}{\ensuremath{\cX}}
\providecommand{\elldual}{{\dual\elldom}}
\renewcommand{\pdt}[1]{{#1}'}

\providecommand{\xespace}{\fespace{\rX}}
\renewcommand{\fes}[1]{\ensuremath{\xespace_{#1}}}

\providecommand{\fesm}[1]{\ensuremath{\xespace_{#1}^\ominus}}
\providecommand{\dgfes}[1]{\ensuremath{\rY_{#1}}}
\renewcommand{\t}[1]{{\ensuremath{t_{#1}}}}
\renewcommand{\tm}[1]{{\ensuremath{t_{#1}^-}}}
\providecommand{\tp}[1]{{\ensuremath{t_{#1}^+}}}
\providecommand{\tpm}[1]{{\ensuremath{t_{#1}^\pm}}}
\providecommand{\tdeg}[1]{{\ensuremath{r_{#1}}}}
\providecommand{\timestep}[1]{{\ensuremath{k_{#1}}}}
\renewcommand{\timestep}[1]{{\ensuremath{\tau_{#1}}}}
\providecommand{\timeinter}[1]{{\ensuremath{I_{#1}}}}
\providecommand{\limposat}[2]{\ensuremath{{#2}(\tp{#1})}}
\providecommand{\limnegat}[2]{\ensuremath{{#2}(\tm{#1})}}
\providecommand{\limpmat}[2]{\ensuremath{{#2}(\tpm{#1})}}
\providecommand{\jumpatof}[2]{\ensuremath{\jump{#2}_{#1}}}
\renewcommand{\ltwop}[3][\pivot]{\qp{#2,#3}_{#1}}
\providecommand{\ltwonorm}[1]{\Norm{#1}_{\pivot}}
\providecommand{\dgenorm}[2][]{\Norm{#2}_{\leb2(\ifx&#1&{\timeinter{n}}\else{#1}\fi;\elldom)}}
\providecommand{\dgdnorm}[2][]{\Norm{#2}_{\leb2(\ifx&#1&{\timeinter{n}}\else{#1}\fi;\dual\elldom)}}
\providecommand{\fe}[1]{{\ensuremath{\mathsf{#1}}}}
\renewcommand{\fe}[1]{{\ensuremath{\mathsf{#1}}}}
\providecommand{\feellop}[1][n]{\ensuremath{\fe A_{#1}}}
\providecommand{\ta}[1]{\ensuremath{\uppercase{#1}}}
\providecommand{\feta}[1]{\fe{\ta{#1}}}
\providecommand{\full}[1]{{\fe{\ta{#1}}}}
\newcommand{\spaceproj}[1][n]{{\fe\pi}_{#1}}

\newcommand{\clementproj}[1][n]{{\fe\Lambda}_{#1}}
\newcommand{\clementprojm}[1][n]{\clementproj[#1]^\ominus}

\newcommand{\fullproj}[1][n]{{\fe\Pi}_{#1}}

\providecommand{\timeosc}[2][]{\ta\Theta\ifx|#1|\else_{#1}\fi{#2}}
\providecommand{\trec}[1]{\widehat{#1}}
\providecommand{\trecta}[1]{\trec{\ta{#1}}}

\providecommand{\trecfeta}[1]{\trec{\feta{#1}}}

\providecommand{\trecerr}{\trec \sigma}
\providecommand{\trecres}{\trec e}

\providecommand{\erec}[1][u]{\widetilde{\uppercase {#1}}}
\providecommand{\erecerr}{\widetilde \sigma}
\providecommand{\erecres}{\widetilde e}
\providecommand{\allres}[1]{\xi_{#1}}
\renewcommand{\abil}[3][]{\duality[#1]{\ellop{#2}}{#3}}
\providecommand{\estfunk}[2][]{\ensuremath{\cE\ifx&#1&{}\else{[#1]}\fi}}
\providecommand{\meshsize}[1][]{\ensuremath{h\ifx|#1|\else_{#1}\fi}}
\renewcommand{\meshsize}[1][]{\ensuremath{h\ifx|#1|\else_{#1}\fi}}
\providecommand{\meshsizem}[1][n]{\ensuremath{h^\ominus\ifx|#1|\else_{#1}\fi}}

\providecommand{\timeind}[1][n]{{\eta_{#1}^{\textup{time}}}}
\providecommand{\spaceind}[1][n]{{\eta_{#1}^{\textup{space}}}}

\providecommand{\oscind}[1][n]{{\eta_{#1}^{\textup{osc}}}}
\providecommand{\dG}{DG\xspace}%
\newboolean{revealconstref}
\setboolean{revealconstref}{false}%
\renewcommand{\constref}[2][]{%
  \ifthenelse{\boolean{revealconstref}}{%
      \ensuremath{\constant{\textup{\ref{#2}{\ifx|#1|{}\else{,\ensuremath{#1}}\fi}}}}
    }{%
      \ensuremath{\constant{\ifx|#1|{}\else{\ensuremath{#1}}\fi}}
}}
\renewcommand{\funkref}[3][]{%
  \ifthenelse{\boolean{revealconstref}}{%
    \ensuremath{{#3}_{\textup{\ref{#2}{\ifx|#1|{}\else{,\ensuremath{#1}}\fi}}}}%
  }{%
    \ensuremath{{#3}_{\textup{{\ifx|#1|{}\else{s\ensuremath{#1}}\fi}}}}%
}}

\providecommand{\memotempa}{}
\providecommand{\memotempb}{}
\usepackage{nicefrac}
\usepackage{bm}
\usepackage{enumitem}
\usepackage{cite}
\usepackage{xcolor}
\hypersetup{
  pdftitle={Parabolic DGCG aposteriori},
  pdfauthor={EH Georgoulis, O Lakkis, TP Wihler},
  pdfsubject={Numerical Analysis},
  pdfkeywords={discontinuous Galerkin, parabolic, partial differential equation, numerical method, finite element, analysis, a posteriori, aposteriori, error estimate, approximation}
}
\usepackage[hyperpageref]{backref}
\provideboolean{showoldtext}
\newcounter{tmpsection}
\newcounter{tmpsubsection}
\newcounter{tmpequation}
\newcommand{\oldtext}[1]{%
  \ifthenelse{\boolean{showoldtext}}{%
    \setcounter{tmpsection}{\value{section}}%
    \setcounter{tmpsubsection}{\value{subsection}}%
    \setcounter{tmpequation}{\value{equation}}%
               {\color{gray!70}%
                 #1}%
    \setcounter{section}{\value{tmpsection}}%
    \setcounter{subsection}{\value{tmpsubsection}}%
    \setcounter{equation}{\value{tmpequation}}%
  }{}%
}
\newcommand{\newtext}[1]{\ifthenelse{\boolean{showoldtext}}{{\color{blue}{#1}}}{#1}}
\newcommand{\tw}[1]{%
  \ifthenelse{\boolean{showchanges}}{{%
      \color{blue}{#1}}}{%
    #1}}%
\renewcommand{\fracl}[2]{\nicefrac{#1}{#2}}
\renewcommand{\jump}[1]{\llbracket #1\rrbracket}

\newcommand{\tnorm}[1]{|\!|\!|#1|\!|\!|}
\renewcommand{\Norm}[1]{\|{#1}\|}
\renewcommand{\memotempa}{%
        \spaceproj%
          \jumpatof{n-1}{%
            f%
          }%
}
\renewcommand{\memotempb}{%
        \jumpatof{n-1}{\spaceproj[]}%
        f%
        (t_{n-1}^+)%
}
\provideboolean{showold}
\provideboolean{showtodo}
\providecommand{\olchanges}[2][]{%
 \ifthenelse{\boolean{showchanges}}{%
   {\ifx|#1|{}\else{%
       \ifthenelse{\boolean{showold}}{%
         \color{black!10!white}{#1}}{%
     }}\fi}%
   {\color{magenta}#2}}{#2}%
}
\providecommand{\twchanges}[2][]{%
	\ifthenelse{\boolean{showchanges}}{%
		{\ifx|#1|{}\else{\color{red}{#1}}\fi}%
		{\color{green!50!red}#2}}{#2}
}
\providecommand{\tand}{,\qquad}
\setboolean{showtodo}{false}%
\setboolean{shownotes}{false}%
\setboolean{showchanges}{false}%
\setboolean{showoldtext}{false}%
\setboolean{showold}{false}%
\setboolean{revealconstref}{false}%

\author{Emmanuil H.~Georgoulis, Omar Lakkis and Thomas P.~Wihler}
\title[A posteriori bounds for fully-discrete $hp$-\dG timestepping methods]{
  A posteriori error bounds for fully-discrete $hp$-discontinuous
  Galerkin timestepping methods for parabolic problems
} 
\address{Emmanuil H.~Georgoulis
  \\
  Department of Mathematics,
  University of Leicester,
  Leicester LE1 7RH,
  United Kingdom
  \\
  and
  \\
  Department of Mathematics,
  School of Applied Mathematical and Physical Sciences,
  National Technical University of Athens,
  Zografou 15780,
  Greece}
\email{\linkedemail{Emmanuil.Georgoulis@le.ac.uk}}
\address{Omar Lakkis
  Department of Mathematics,
  University of Sussex,
  Falmer,
  Brighton,
  GB-BN1 9QH,
  England UK}
\email{\linkedemail{lakkis.o.maths@gmail.com}}
\address{
  Thomas~P.~Wihler
  \\
  Mathematisches Institut,
  Universit\"at Bern,
  Sidlerstrasse 5,
  CH-3012 Bern,
  Switzerland}
\email{\linkedemail{wihler@math.unibe.ch}}
\date{\today}
\begin{document}
\begin{abstract}
  We consider fully discrete time-space approximations of
  abstract linear parabolic partial differential equations (PDEs)
  consisting of an $hp$-version discontinuous Galerkin (\dG) time
  stepping scheme in conjunction with standard (conforming) Galerkin
  discretizations in space. We derive abstract computable
  \aposteriori error bounds resulting, for instance, in concrete
  bounds in $\leb{\infty}(I;\leb2(\W))$- and
  $\leb{2}(I;\sobh{1}(\W))$-type norms when $I$ is the temporal and
  $\W$ the spatial domain for the PDE. We base our methodology for
  the analysis on a novel space-time reconstruction approach.  Our
  approach is flexible as it works for any type of elliptic error
  estimator and leaves their choice to the user. It also
  exhibits mesh-change estimators in a clear and concise way.
  We also show how our approach allows the derivation of such bounds
  in the $\sobh1(I;\sobh{-1}(\W))$ norm.
\end{abstract}
\maketitle
\section{Introduction}
Adaptive numerical methods have been shown to
provide accurate and efficient
numerical treatment of evolution PDEs
thanks to their properties for localized mesh resolution especially in the
context of moving fronts, interfaces, singularities, or layers (both
boundary and interior). Such numerical methods predominantly admit
spatial discretizations of variational type, e.g., finite element
methods (FEMs), which allow for general,
possibly unstructured, dynamic mesh modification. FEMs are also
ideally suited for deriving mathematically rigorous \aposteriori
bounds, owing to their variational nature.

Some of the classical works on adaptive finite
element methods for parabolic problems
\cite{EJ-I,EJ-II,EJ-IV,EJ-V,EJ-VI} are based on discontinuous Galerkin
(\dG) time stepping combined with FEM in space, and proving
\aposteriori bounds in various norms using duality techniques. The key
motivation in using \dG in time, which is also of variational type, is
that it naturally allows for spatially-local-time stepping, i.e.
different time step sizes in different parts of the spatial domain
\cite{J78,EJT85,EJ-I, MB:97}. This classical, but as of yet
undeveloped in full, concept of local adaptivity in both space and
time has the potential of delivering substantial computational savings
and even complexity reduction.

In addition to the ability of Galerkin time marching schemes to
employ locally different time step sizes, their variational
character also allows for arbitrary variations in the local
approximation orders. They can therefore be cast naturally into the
framework of $hp$-approximation schemes. In the context of
parabolic PDEs, $hp$-version time marching methods can be
used, for instance, to resolve an initial layer in the (otherwise
smooth) solution at high algebraic or even exponential rates of
convergence, see, e.g., the
works~\cite{SchoetzauSchwab00,SchoetzauSchwab01,Gerdes} on linear
parabolic PDEs, and also~\cite{MSW1,MSW2} which employ a
combination of~$hp$-version time stepping with suitable wavelet
spatial discretizations to yield a log-linear complexity algorithm
for nonlocal evolution processes involving pseudo-differential
operators. Additionally, we note the numerical analysis of
high-dimensional parabolic problems using sparse grids in space;
see~\cite{SchwabPetersdorff}.

More recent results on rigorous \aposteriori bounds for parabolic
problems have focused on extending the paradigm of the reliable and
efficient \aposteriori error analysis of elliptic problems to the
parabolic case~\cite{P98,V98,V03}. Such works typically involve basic
low-order time stepping schemes combined with various types of FEM in
space. \Aposteriori error bounds for \dG time-stepping methods have
also appeared in the last few years; we point
to~\cite{MakNoc06,SW10,KMW:16} which are based on the
\emph{reconstruction} technique, to~\cite{ESV:16,ESV:17} which employ
an equilibrated flux approach, or to~\cite{GKSZ17} which presents a
provably convergent adaptive algorithm for a residual-type a
posteriori estimator.

In this paper, we present \aposteriori
error bounds for an $hp$-version \dG-in-time and conforming Galerkin
discretization in space method for \emph{both}
  $\leb{\infty}(I;\leb2)$- and $\leb{2}(I;\sobh{1})$-norm errors
  separately, allowing for what appears to be optimal order in each
  case.  The key idea is the use of suitable reconstruction frameworks
  to derive a perturbed PDE for the reconstructed error of the numerical method;
\aposteriori error bounds are then deduced using PDE stability properties,
cf.~\cite{MakNoc03,MakNoc06,LakMak06, AMN06}. Our approach is based
on new space and space-time reconstructions which are built
  on the combination of respective ideas for \dG-time stepping methods
\cite{MakNoc06,SW10} and elliptic reconstruction \cite{MakNoc03,
  LakMak06} to the fully-discrete setting. To that end, the key
challenge of constructing a globally time-continuous reconstruction in
the presence of mesh modification between time-steps is addressed by
first reconstructing onto the solution space with respect to the
spatial variables via a novel elliptic reconstruction
  definition, given in
  \eqref{eqn:def:elliptic-reconstruction:strong}, which is a
  modification of the one in \cite{LakMak06}. In particular, the new proposed
  elliptic reconstruction takes into account the effect of \emph{mesh-change}. 
  
Our results are closely related, however, with important departures,
to those of \cite{ESV:16,ESV:17, GKSZ17}.  In particular, the new reconstructions defined below allow for the derivation
of \aposteriori upper error bounds for each of the following norms
$\leb2(I;\elldom)$ (Theorem~\ref{thm:L2H1}) and
$\sobh1(I;\dual\elldom)$ (\S \ref{sec:H1(Hm1)apost})
separately; the Hilbert space $\elldom$ is the domain of a
self-adjoint uniformly elliptic operator $\ellop$ (see~\S\ref{sec:setup}
for details). A key attribute
of our approach is the flexibility in incorporating any \aposteriori
elliptic error estimators available, as the reconstruction-type approach, allows to separate the challenges in the \aposteriori error estimation of elliptic and time-evolution errors. To facilitate a wide range of applications, we will present the theory within a Gelfand-triple-type abstract setting allowing, for instance, both second- and fourth-order spatial operators. This generality comes at the possible expense of different, yet quantitatively analogous, computable constants in the resulting \aposteriori error estimators compared to the bounds in \cite{ESV:16,ESV:17, GKSZ17}. For instance, when the
equilibrated \tw{flux} elliptic error estimators from \cite{BPS:09}
are used (with $\elldom=\sobhz1(\W)$ and $\pivot=\leb2(\W)$), we
recover similar estimators to the upper bounds derived in
\cite{ESV:16,ESV:17}. Importantly, however, the work \cite{ESV:16}
shows that these are also lower bounds for the ``joint-norm'' of
$\sobh1(I;\dual\elldom)\meet\leb2(I;\elldom)$, and the article
\cite{ESV:17} does the same for the $\leb2(I;\elldom)$ under the
condition $h^2<c\tau$, relating the mesh-size $h$ with the time-step
$\tau$ for some constant $c>0$. Also, in the present work, we are
not concerned with the interesting question of convergence of
adaptive algorithms as in \cite{GKSZ17}. Crucially, however, the novel space-time reconstruction, allows for the proof of an \aposteriori error bound for the
$\leb\infty(I;\pivot)$-norm, which appears to be of optimal order; this result, to the best of our reading, is not captured in \cite{ESV:16,ESV:17,GKSZ17}.
\subsubsection*{Outline}
The remainder of this work is structured as
follows. In \S\ref{sec:setup} we set up the abstract framework for the paper by
introducing the model parabolic PDE problem and its \dG-in-time and conforming Galerkin spatial discretization. Furthermore, in
\S\ref{sec:reconstructions}, we provide the necessary technical tools for the ensuing analysis, and state their essential
properties. In \S\ref{sec:L2(H1)apost}, we derive
  \aposteriori error bounds in the $\leb{2}(I;\elldom)$-norm using a
  time reconstruction and a novel
  elliptic reconstruction which includes mesh-change; this technical novelty
  is revealed to be crucial in our setting.
In \S\ref{sec:Linfty(L2)apost}, upon defining a new \emph{space-time} reconstruction, we derive \aposteriori
  error bounds in the $\leb\infty(I;\pivot)$-norm. Finally, in
  \S\ref{sec:H1(Hm1)apost} we briefly discuss an approach to arriving
  at $\sobh1(I;\dual\elldom)$-type \aposteriori error estimates.
\section{Model problem and space-time discretization}
\label{sec:setup}
We introduce most of the notation and technical background for the
paper.  In~\S\ref{sec:abstract-setting} we provide the functional
analytic set-up for the abstract heat equation, a related concrete
Example~\ref{exa:concrete-elliptic-operators}, and we present the fully discrete numerical scheme in \S\ref{sec:Galerkin-scheme}.
\subsection{Abstract setting}
\label{sec:abstract-setting}
Throughout this work, Bochner spaces will be used. To that end, given
an interval $J\subset\R{}$, and a real Hilbert
space~$\cZ$ with inner product $\ltwop[\cZ]\cdot\cdot$ and induced
norm~$\Norm\cdot_{\cZ}$, we define
\begin{equation}
  \Norm{u}_{\leb{p}(J;\cZ)}
  =
  \begin{cases}
    \displaystyle\left(\int_J\Norm{u(t)}^p_{\cZ}\d t\right)^{\nicefrac{1}{p}},
    &
    1\le p<\infty,
    \\
    \esssup_{t\in J}\Norm{u(t)}_{\cZ},&p=\infty.
  \end{cases}
\end{equation}
We write $\leb{p}(J;\cZ)$ to signify the space of measurable
functions $u:J\to\cZ$ such that the corresponding norm is
bounded. Note that~$\leb2(J;\cZ)$ is a Hilbert space with inner
product and induced norm given by
\begin{gather}
  \ltwop[{\leb2(J;\cZ)}]uv=\int_J(u(t),v(t))_{\cZ}\d t,
\end{gather}
and $ \Norm{u}_{\leb2(J;\cZ)}:= \ltwop[{\leb2(J;\cZ)}]{u}{u}^{\nicefrac12}$,
respectively. We also let~$\sobh1(J;\cZ)$ be the Sobolev space of
all functions in~$\leb2(J;\cZ)$ whose (weak temporal) derivative is
bounded in~$\leb2(J;\cZ)$, with the norm
\begin{equation}
  \Norm{u}_{\sobh1(J;\cZ)}
  =
  \powqp{\nicefrac12}{
    \int_J\left[
    \Norm{u(t)}_{\cZ}^2+\Norm{u'(t)}_{\cZ}^2\right]
    \d t}.
\end{equation} 
Finally, the space $\cont0(\closure{J};V)$ consists of all
functions that are continuous on~$\closure{J}$, the closure of $J$,
with values in
$\cZ$, endowed with the standard maximum norm
\begin{equation}
  \Norm{u}_{\cont0(\closure{J};\cZ)}
  =
  \max_{t\in \closure{J}} \Norm{u(t)}_{\cZ}.
\end{equation}
We consider henceforth two (real) Hilbert spaces $\elldom$ and
$\pivot$ forming a Gelfand triple
\begin{equation}\label{eq:Gelfand}
  \elldom\embedsin \pivot\embedsin\dual \elldom,
\end{equation}
where $\dual\elldom$ denotes the dual of $\elldom$. The duality
pairing $\duality\cdot\cdot$ of~$\dual\elldom$ and~$\elldom$ can be
seen as a continuous extension of the inner
product~$\ltwop{\cdot}{\cdot}$. In particular,
identifying~$\dual\pivot\simeq\pivot$, for~$u\in\pivot$
and~$v\in\elldom$, there holds
\begin{equation}\label{eq:identify}
  \duality{u}{v}=\ltwop{u}{v};
\end{equation}
see, e.g.,~\cite[\S7.2]{Roubicek:13} for details.

Moreover, let
\begin{equation}\label{eq:A}
  \funk\ellop \elldom{\dual \elldom}
\end{equation}
be a \emph{self-adjoint linear elliptic operator}
continuous and coercive
in the sense that there exist constants
$\beta\geq\alpha>0$
such that
\begin{equation}
\label{eq:ellop-prop}
\begin{split}
\duality{\ellop v}w&\leq
\beta
\Norm v_\elldom\Norm w_\elldom\Foreach v,w\in\elldom,
\\
\duality{\ellop v}v
&
\geq
\alpha
\Norm v_\elldom^2
\Foreach v\in\elldom.
\end{split}
\end{equation}

Given an initial value $u_0\in \pivot$, a final time~$T>0$,
denoting henceforth the time interval
\begin{equation}\label{eq:I}
  I:=\opclinter0T,
\end{equation}
and given a source function $f\in\leb2(I;\elldual)$,
we are interested in a
Galerkin-type numerical approximation of the function
\begin{equation}\label{eq:model}
  u\in\sobh1(I;\dual\elldom)\meet\leb2(I;\elldom),
\end{equation}
which solves uniquely the linear parabolic initial value problem
\begin{equation}
  \label{eqn:def:exact-pde}
  \pdt u+\ellop u=f
  \AND
  u(0)=u_0.
\end{equation}
Incidentally, due to the continuous embedding
\begin{equation}
  \sobh1(I;\dual\elldom)\meet\leb2(I;\elldom)
  \embedsin
  \cont0(\closure I;\pivot),
\end{equation}
it follows that~$u$ belongs to $\cont0(\closure I;\pivot)$
\cite[e.g.,~Lemma~7.3]{Roubicek:13} and the initial condition
in~\eqref{eqn:def:exact-pde} makes sense.
\begin{Exa}[concrete elliptic operators]
  \label{exa:concrete-elliptic-operators}
  A commonly encountered situation which can be cast in the above
  framework is the classical linear diffusion equation, i.e. $\ellop
  v=-\nabla\cdot[\mat A\grad v]$, where, for a given open, connected, and
  bounded domain $\W\subset\R d$, $d=1,2$ or $3$, we consider a given
  symmetric matrix-valued function~$\funk{\mat A}\Omega{\realmats dd}$,
  $\mat{A}\in\leb\infty(\Omega)^{d\times d}$, satisfying
  \begin{equation}
    \label{const:concrete-coercivity}
    \transposevec v\mat A(\vec x)\vec v
    \ge
    \alpha
    \norm{\vec v}^2
    \Foreach
    \vec x
    \in
    \Omega
    ,
    \Foreach
    \vec v
    \in
    \R{d},
  \end{equation}
  for some constant \(\alpha>0\). Here, we choose, e.g.,
  $\pivot:=\leb2(\W)$, $\elldom:=\honezw$, and
  $\elldual:=\sobh{-1}(\W)$ to be the typical function spaces in the
  context of second-order linear elliptic PDEs with homogeneous Dirichlet boundary condition $u=0$ on $\boundary\W$. 
 
Another possible choice, e.g., is  $\pivot:=\leb2(\W)$, $\elldom:=\sobh{2}_0(\W)$, and
    $\elldual:=\sobh{-2}(\W)$, for the case of a fourth-order parabolic problem with essential boundary conditions. 
\end{Exa}
\subsection{Time discontinuous and space conforming Galerkin approximation}
\label{sec:Galerkin-scheme}   
Given a (real) linear space $\cZ$, the space of all $\cZ$-valued polynomials of degree at most~$r$, with~$r\in\mathbb{N}_0$, on~$\reals$ is defined by
\begin{equation}
  \begin{gathered}
    \poly r(\cZ)
    :=
    \textstyle
    \ensemble{\funk p{\reals}{\cZ}}{p(x)
      =
      \sum_{i=0}^r z_{i}x^i\text{ for some }(z_0,\ldots,z_r)\in\cZ^{r+1}}.
  \end{gathered}
\end{equation}
In addition, if $D\subseteq\reals$, we
define
\begin{equation}
  \poly r(D;\cZ):=\ensemble{\restriction pD}{p\in\poly r(\cZ)}.
  \end{equation}

In order to introduce the discontinuous Galerkin time stepping scheme
for~\eqref{eqn:def:exact-pde}, we consider a finite sequence of
\emph{time nodes} and \emph{time steps},
\begin{equation}
  0=\t0<\t1<\dotsc<\t N=T,    
  \AND
  \timestep n:=\t n-\t{n-1}\text{ for }\rangefromto n1N,
\end{equation}
as well as the corresponding \emph{time intervals}
\begin{equation}\label{eq:cI}
  \timeinter n:=
  \begin{cases}
    \setof 0
    &
    \text{ for }n=0,
    \\
    \opclinter{\t{n-1}}{\t n}
    &
    \text{ for }\rangefromto n1N.
  \end{cases}
\end{equation}
Thus, we have a \emph{partition}
$\cI:=\ensemble{\timeinter n}{\rangefromto{n}{1}{N}}$
of the time interval $I$ given in~\eqref{eq:I}.

Given a $\cI$-piecewise continuous function $\funk g{I\subseteq\reals}{\cZ}$,
we define its time jump across $t_n$, $\rangefromto n0{N-1}$, for given~$g(t_0^-)$, by
\begin{equation}\label{eq:jump}
  \jumpatof {n}g
  :=
  \limposat ng
  -
  \limnegat ng,
\end{equation}
where we introduce the one-sided limits
$
\limpmat{n}g
  =
  \lim_{\epsi\to0^+}g(t_n\pm\epsi).
$
Moreover, we associate with the finite sequence of time instants
$\listidotsfromto t0N$ a finite sequence of \emph{finite-dimensional
  conforming subspaces}
\begin{equation}
  \label{eq:def:finite-dimensional-conforming-subspaces}
  \fes n\subset\elldom,\text{ for }\rangefromto n0N.
\end{equation}

For a generic $\elldom$-conforming Galerkin space $\xespace$, we
signify by $\spaceproj[\xespace]$ the $\pivot$-\emph{orthogonal
  projection} from $\dual\elldom$ onto~$\xespace$:
\begin{equation}
  \label{eqn:def:l2projection-from-pivot}
  \dfunkmapsto[.]{
    \spaceproj[\xespace]
  }{
    v
  }{
    \dual\elldom
  }{
    \spaceproj[\xespace] v
    :\ \ltwop{\spaceproj[\xespace] v}{\fe w}=\duality{v}{w}
    \Foreach\fe w\in\xespace
  }{
    \xespace
  }
\end{equation}
Note that, due to~\eqref{eq:identify}, for~$v\in\pivot$, we have
\begin{equation}
  \ltwop{\spaceproj[\xespace] v}{\fe w}
  =\duality{v}{\fe w}=\ltwop{v}{\fe w}\Foreach\fe w\in\xespace.
\end{equation}
When $\xespace$ is $\fes n$, for some
$n\integerbetween0N$, we write $\spaceproj$ to indicate $\spaceproj[\xespace]$.

In order to introduce the time semidiscrete and space-time fully
discrete spaces, let $r_n\in\NO$, $\rangefromto n1N$, be a polynomial
degree. Then, consider the 
\emph{time semidiscrete Galerkin space}
\begin{equation}
  \cY
  :=
  \ensemble{
    \funk{\full V}{\opclinter0T}{\elldom}
  }{
    \restriction{\full V}{\timeinter n}
    \in
    \poly{\tdeg n}(\timeinter n;\elldom)
    \Foreach\rangefromto n1N
},
\end{equation} 
respectively, 
the \emph{space-time fully discrete Galerkin space}
\begin{equation}
  \dgfes{} 
:=
\ensemble{
  \funk{\full V}{\opclinter0T}{\elldom}
}{
  \restriction{\full V}{\timeinter n}
  \in
  \poly{\tdeg n}(\timeinter n;\fes n)
  \Foreach\rangefromto n1N
}
,
\end{equation} 
where $\poly{\tdeg n}(\timeinter n;\fes n)$ are the space-time
Galerkin subspaces. The \emph{fully discrete time-discontinuous
  Galerkin and spatially-conforming} approximation
of~\eqref{eqn:def:exact-pde} is then an $\cI$-piecewise continuous function
$\full u\in\dgfes{}$, 
such that
\begin{gather}
  \full u(\tm 0):=\spaceproj[0] u_0,\label{eq:dG1}
\end{gather} 
and for $\rangefromto n1N$,
\begin{equation}
  \label{eq:dG2}
  \begin{split}
    \int_\timeinter n
    \qb{
      \ltwop{\pdt{\full u}}{\full v}
      +
      \abil{\full u}{\full v}
    }\d t
    &+
    \ltwop{\jumpatof {n-1}{\full u}}{\full v(\tp {n-1})}
    =
    \int_\timeinter n
    \duality{f}{\full v}\d t,
  \end{split}
\end{equation}
$ \Foreach
\full v\in\poly{\tdeg n}(\timeinter n;\fes n),$ where~$\jumpatof{0}{\full U}=\full U(t_0^+)-\spaceproj[0] u_0$.
\section{Reconstructions}
\label{sec:reconstructions}
We will next introduce some technical essentials.  The
main tools are the time
lifting (\S\ref{def:time-lifting}), the time reconstruction
(\S\ref{def:time-reconstruction}), and a \emph{new variant} of the elliptic reconstruction from \cite{LakMak06} for fully discrete schemes (\S\ref{sec:ellrec}).  In
\S\ref{the:elliptic-reconstruction-error-estimates} we postulate the availability of \aposteriori error estimators for elliptic residuals, and we give some
pointers to the relevant literature.  In addition, we discuss various
error estimates that measure the time reconstruction error; in
particular, we state two identities which follow directly,
respectively, from~\cite[Theorem~2]{SW10} and, taking into account the
explicit representation of the time reconstruction,
from~\cite[Lemma~1]{HolmWihler:15}.
\subsection{Time lifting}
\label{def:time-lifting}
Let us consider, for given $\rangefromto{n}{1}{N}$, a linear \emph{time lifting operator} 
\begin{equation}\label{eq:tl1}
\chi_n:
\pivot
\to 
\poly{\tdeg n}(\timeinter n;\pivot).
\end{equation}
It is defined, for each $w\in \pivot$, by the Riesz representation
\begin{equation}\label{eq:tl2}
\int_\timeinter n    
\ltwop{\chi_n(w)}{\ta v}
=
\ltwop{w}{\ta v(\tp {n-1})} \quad 
\Foreach
\ta v\in\poly{\tdeg n}(\timeinter n;\pivot).
\end{equation}
\begin{Lem}[space invariance under time lifting]
	For any linear subspace $\cW\subseteq\pivot$, the time lifting from~\eqref{eq:tl1} and~\eqref{eq:tl2} satisfies
	\begin{equation}
	w\in\cW\implies\chi_n(w)\in\poly{\tdeg n}(I_n;\cW)
	.
	\end{equation}
	In particular, writing $\charfun A$ for the indicator function
	on a generic set $A$, and assuming $\fe w_n\in\fes n$ for
        each $n\integerbetween1N$,
	we have 
	\begin{equation}\label{eq:indicatorfct}
	\sumifromto n1N\chi_n(\fe w_n)\charfun{\timeinter n}
	\in\dgfes.
	\end{equation}
\end{Lem}
\begin{proof}
	This result is a straightforward consequence of the explicit
	representation of~$\chi_n$ as described in~\cite[Lemma~6]{SW10}.
\end{proof}
\subsection{Time reconstruction}
\label{def:time-reconstruction}
Let us define the \emph{time-reconstruction} $\trecta w$
of a given time-discrete function 
\begin{equation}
\ta w \in
  \ensemble{
    \funk{\full V}{\opclinter0T}{\pivot}
  }{
    \restriction{\full V}{\timeinter n}
    \in
    \poly{\tdeg n}(\timeinter n;\pivot)
    \Foreach\rangefromto n1N
},
\end{equation} 
as follows: for
each $n\integerbetween1N$, we let
$\trecta w|_{\timeinter n}\in\poly{\tdeg{n}+1}(\timeinter n;\pivot)$ satisfy
\begin{equation}
\label{eq:time-reconstruction-formula}
\restriction{\trecta w}{\timeinter n}(t)
:=
\limnegat{n-1}{\ta w}
+
\int_{t_{n-1}}^t\qbbig{\pdt{\ta w}(s)+\chi_n(\jumpatof{n-1}{\ta w})(s)}\d s
\qquad\text{ for }
t\in\closure{I}_n.
\end{equation}
Equivalently, we note the following characterization of $\trecta w$ in weak form
on each $\timeinter n$,
\begin{gather}
\label{eq:time-reconstruction-definition}
\begin{split}
\int_\timeinter n(\pdt{\trecta w},\ta v)_{\pivot}\d t
=
\ltwop{\jumpatof {n-1}{\ta w}}{\limposat{n-1}{\ta v}}
+
\int_\timeinter n\ltwop{\pdt{\ta w}}{\ta v}\d t,
\\
\end{split}
\intertext{$\Foreach\ta  v\in \poly{\tdeg n}(\timeinter n; \pivot)$, with the initial condition}
\limposat {n-1}{\trecta w}
:=
\limnegat{n-1}{\ta w},
\end{gather}
for $n\integerbetween1N$; cf.~\cite{MakNoc06, SW10}. Evidently, the above construction carries over to any linear subspace $\cW\subset\pivot$ in an obvious way.
\begin{Pro}[time-reconstruction error identities]
	\label{pro:time-reconstruction-error-estimate}
	Consider any (real) Hilbert space $\cW$, and
	$\restriction{\ta w}{\timeinter n}\in\poly{\tdeg{n}}(\timeinter{n};\cW)$, 
	$\rangefromto{n}1N$, with $\trecta{w}$ defined from $\ta{w}$ through
	\eqref{eq:time-reconstruction-formula}. Then, for given~$\limnegat{0}{\ta w}\in\cW$, the following approximation
	identities hold
	\begin{equation}\label{eq:rec1st}
	\Norm{\ta W -\trecta w}_{\leb2(\timeinter n;\cW)}
	=
	\constref[\timestep n,\tdeg n]{eqn:time-reconstruction:constant}
	\Normreg{\jumpatof{n-1}{\ta w}}_{\cW},
	\end{equation}
	where
	\begin{equation}
	\label{eqn:time-reconstruction:constant}
	\constref[\timestep{},\tdeg{}]{eqn:time-reconstruction:constant}
	:=
	\left(\frac{\timestep{}(\tdeg{}+1)}{(2\tdeg{}+1)(2\tdeg{}+3)}\right)^{\nicefrac12},
	\end{equation}
	and
	\begin{equation}\label{eq:time-reconstruction-linf}
	\Norm{\ta w -\trecta w}_{\leb{\infty}(\timeinter n;\cW)}
	=\Normreg{\jumpatof{n-1}{\ta w}}_{\cW}.
	\end{equation}
\end{Pro}
\begin{proof}
	The identity~\eqref{eq:rec1st} was first proven in~\cite[Lemma 2.2]{MakNoc06},
	and extended to this exact form in~\cite[Theorem~2]{SW10} accounting
	for the dependence on the polynomial degree explicitly. The second
	equality~\eqref{eq:time-reconstruction-linf} follows directly by combining the explicit representation
	formula derived in~\cite[Eq.~(33)]{SW10}
	with~\cite[Lemma~1]{HolmWihler:15}.
\end{proof}

	\begin{Obs}[continuity of the time-reconstruction]
		\label{obs:discontinuity-of-time-reconstruction}
		Owing to \cite[Lemma~2.1]{MakridakisNochetto:06:article:A-posteriori}, the semidiscrete (spatially exact) time-reconstruction \eqref{eq:time-reconstruction-formula} originally defined
		in~\cite{MakridakisNochetto:06:article:A-posteriori} and \cite{SW10} is a continuous function in time. In particular, the time-reconstruction $\trecfeta u$ of the fully discrete solution $\full u$, defined in~\eqref{eq:dG1} and~\eqref{eq:dG2}, is still continuous across the time nodes $t_0,\ldots,t_{N-1}$, despite having $\spaceproj\trecfeta{U}\neq \trecfeta{U}$, on $I_n$, when the spatial mesh changes across $t_{n-1}$ in a non-hierarchical fashion.
	\end{Obs}
\subsection{Elliptic reconstruction}
\label{sec:ellrec}
Let $\xespace\subset\elldom$ be a generic conforming Galerkin
space. Then, given the elliptic operator~$\ellop$ from~\eqref{eq:A},
we define the \emph{discrete elliptic operator}
$\funk{\feellop[\xespace]}{\xespace}{\xespace}$, for each $\fe
w\in\xespace$, as \(\feellop[\xespace] \fe w\in\xespace\) such that
\begin{equation}\label{eq:Adiscrete}
   \ltwop{\feellop[\xespace]\fe w}{\fe v}
  =
  \duality{\ellop\fe w}{\fe v}
  \Foreach
  \fe v\in\xespace.
\end{equation}
From the ellipticity of $\ellop$, it follows that $\funk{\feellop[\xespace]}{\xespace}{\xespace}$ is
invertible. 
Note that the discrete elliptic operator's domain may be extended
from $\xespace$ to all of~$\elldom$; indeed, this may be convenient in some cases where
we are ready to give up its invertibility. If $\xespace$ is one of $\fes n$s, for
some~$\rangefromto{n}0N$, we denote $\feellop[\fes n]$
by $\feellop$. 

 To optimize on the structure of the mesh-change indicator
  below, we use a nonstandard elliptic reconstruction on each time
  interval $I_n$. To that end, for each $t\in I_n$, we define the
  \emph{elliptic reconstruction} $\erec(t)\in \elldom$, by
  \begin{equation}\label{eqn:ell_rec_mod}
    \duality{\ellop\erec(t)}v
    =
    \ltwop{\feellop\full u(t)+ \spaceproj\trecfeta{U}'(t)-\trecfeta{U}'(t)}v
    \quad  \Foreach
    v\in\elldom,
  \end{equation}
  where $\trecfeta u$ is the time-reconstruction of the fully discrete solution $\full u$ from~\eqref{eq:dG1} and~\eqref{eq:dG2}.
  It follows that $\erec\in\cY$ and may be written implicitly, for
  any~$\rangefromto{n}{1}{N}$, as the solution of the $t$-dependent
  elliptic problem
  \begin{equation}
    \label{eqn:def:elliptic-reconstruction:strong}
    \ellop\erec(t)
    =
    \feellop 
    \full u (t)+ \spaceproj\trecfeta{U}'(t)-\trecfeta{U}'(t)
    \quad  \text{ for }
    t\in I_n.
  \end{equation}
The initial value of $\erec$ is given by
\begin{equation}
  \erec(0)=\limnegat0{\erec}:=
  \spaceproj[0]u(\tm{0})=\spaceproj[0]u_0,
\end{equation}
with~$u_0\in\pivot$ from~\eqref{eqn:def:exact-pde}.

Upon restricting the test functions in \eqref{eqn:ell_rec_mod} to $\fes n$, we obtain
\begin{equation}\label{eqn:erec_orthog}
  \abil{\erec(t)}{\fe v}
  =
  \ltwop{\feellop\full u(t)}{\fe v}
  =
  \abil{\full u(t)}{\fe v}
\quad  \Foreach
  \fe v\in\fes n,
\end{equation}
cf.~\eqref{eq:Adiscrete}. Evidently this identity implies
\begin{equation}
  \label{eqn:elliptic-reconstruction-identity:space-time}
  \int_{\timeinter n}
  \abil{\erec}{\full v}\d t
  =
  \int_{\timeinter n}
  \ltwop{\feellop\full u}{\full v}\d t
 ,
\end{equation}
for each $\full v\in\poly{\tdeg n}(\timeinter n;\fes n)$, $\rangefromto{n}{1}{N}$.

\subsection{Assumption (elliptic \aposteriori error estimates)} 
\label{the:elliptic-reconstruction-error-estimates}
For given~$g\in\pivot$, consider the abstract elliptic problem of finding $w\in\elldom$ such that $\ellop w=g$. Moreover, let $\xespace\subseteq\elldom$ be a generic $\elldom$-conforming Galerkin space, and let $\fe w\in\xespace$
be $w$'s Galerkin approximation in $\xespace$, defined implicitly as
the solution of $\feellop[\xespace] \fe w=\spaceproj[\xespace] g$. Then some \emph{\aposteriori error bound} holds, viz.,
\begin{equation}
  \label{eq:ellapost}
  \Norm{w -\fe w}_{\cZ}
  \leq
  \cE_{\cZ,\xespace}[\fe w, g],
\end{equation}
with a suitable \emph{\aposteriori error estimator} $\cE_{\cZ,\xespace}$, which
we assume to be available for $\cZ$ representing any of the spaces
$\elldom,\pivot$ or $\dual\elldom$.
Recalling~\eqref{eqn:def:elliptic-reconstruction:strong}, assumption
\eqref{eq:ellapost} allows, for instance, to get \aposteriori error
control of the~\emph{elliptic reconstruction error} in the $\cZ$-norm,
i.e. 
\begin{equation}\label{eq:elliptic_apost}
  \|\erec  -\full u\|_{\cZ}
  \leq
  \cE_{\cZ,\xespace_n}[\full u, \feellop \full u +\spaceproj\trecfeta{U}'-\trecfeta{U}']\quad\text{on }I_n,
\end{equation}
for the selection of spaces~$\cZ=\elldom,\pivot,\text{or
}\dual\elldom$.  Details on such \aposteriori error estimates can be
found, e.g., in \cite{AinsworthOden:00:book:A-posteriori,Braess:07:book:Finite,BrennerScott:07:book:The-mathematical,DM:99}.  It is worth mentioning
that~$\fe w$ does not need to belong to $\xespace$
for~\eqref{eq:ellapost} to hold; instead, it is usually enough that
$\qp{w-\fe w}$ is $\ellop$-orthogonal to $\xespace$ in order to derive
elliptic \aposteriori error estimates.
  \subsection{Pointwise form}
  For~$\rangefromto{n}{1}{N}$, we denote by
  \begin{equation}
  \fullproj:\,\leb 2(\timeinter n;\elldom')\to\poly{\tdeg n}(\timeinter n;\fes
  n),\qquad f\mapsto\fullproj f,
  \end{equation}
   the time-local fully discrete~$\leb{2}(\timeinter n;\pivot)$-orthogonal
    projection defined by
    \begin{equation}
      \int_{\timeinter n}\ltwop{\fullproj f}{\full v}\d t
      =
      \int_{\timeinter n}\duality{f}{\full v}\d t\
      \Foreach \full v\in\poly{\tdeg n}(\timeinter n;\fes n),
    \end{equation}
    and its time-global counterpart
    \index{time-global fully discrete~$\leb2(I;\pivor)$-orthogonal projection}
    \begin{equation}
      \dfunkmapsto[.]{
        \fullproj[]
      }
      f{
        \leb 2(I;\elldom')
      }{
        \fullproj[]f:\restriction{\fullproj[]f}{I_n}=\fullproj f\Foreach n=1,\dots,N
      }{
        \dgfes{}
      }
    \end{equation}
  
  Let $\full V\in\dgfes{}$, and note that using
  Definition \ref{def:time-reconstruction}
  and identity \eqref{eqn:elliptic-reconstruction-identity:space-time},
  the fully discrete local DG formulation~\eqref{eq:dG2} transforms into
  \begin{equation}
    \int_{\timeinter n}
    \ltwop{
      \pdt{\trecfeta u}
      +
      \feellop\full u
      -
      f
    }{\full v}\d t
    =
    0,\qquad n=1,\ldots,N.
  \end{equation}
  Thus,  for each $\rangefromto n1N$, we have
  \begin{equation}
    \label{eqn:time-reconstruction-pde:pointwise-form}
    \spaceproj\pdt{\trecfeta u}( t)
    +
    \feellop \full u( t)
    -
    \fullproj f(t)
    =
    0\quad
    \Foreach t\in\timeinter n,
  \end{equation}
  noting that $\fullproj\pdt{\trecfeta u}=\spaceproj\pdt{\trecfeta u}$
  from Fubini's Theorem. Hence, from the elliptic
  reconstruction \eqref{eqn:ell_rec_mod}, we deduce the pointwise form
  \begin{equation}
    \label{eqn:reconstruction-pdr:pointwise-form}
    \pdt{\trecfeta u}(t)
    +
    \ellop \erec(t)=
    \fullproj  f(t)
    \ \text{ $\ $for each }
    t\in\timeinter n
    .
  \end{equation}
  \begin{Obs}\label{rem:magic_identity}
    Upon noting the trivial identity for $t\in I_n$,
    \begin{equation}
      \ellop \erec(t)=
      \fullproj  f(t)-\pdt{\trecfeta u}(t),
    \end{equation}	
    we observe that the above definition of the elliptic reconstruction constitutes a high order extension of the zero-th order version from \cite[Definition 6.1]{CGM}.
  \end{Obs}
  \section{$\leb2(I;\elldom)$-norm \aposteriori error analysis}\label{sec:L2(H1)apost}
  To highlight the potential generality of the reconstruction
  approach, we first embark on proving \aposteriori error bounds in
  the $\leb2(I;\elldom)$-norm. As we will see below, the resulting
  bounds are qualitatively closely related to the bounds from
  \cite{ESV:16,ESV:17,GKSZ17} when $\elldom = \sobh{1}(\W)$. A key attribute of the approach
  presented below is that the resulting \aposteriori bounds are
  flexible with respect to the choice of respective elliptic bounds,
  such as residual, recovery, or flux-reconstruction ones.
  \subsection{Errors}
  Introduce the following \emph{errors}
  \begin{equation}
    \trecerr
    :=
    \trecfeta u-\full u
    \tand
    \erecerr
    :=
    \erec -\full u
  \end{equation}
  and the corresponding \emph{remainder error}
  \begin{equation}
    \trecres
    :=
    u-\trecfeta u
\    \AND \ 
    \erecres
    :=
    u-\erec
    ,
  \end{equation}
  whereby the \emph{full error} can be decomposed as
  \begin{equation}
    e
    :=
    u-\full u
    =
   \trecres +\trecerr
    =
    \erecerr+\erecres.
  \end{equation}
  Recalling that our ultimate goal is to find \aposteriori estimates
  for $e$, and noting that such estimates for $\trecerr$ and
  $\erecerr$ are provided by Proposition
  \ref{pro:time-reconstruction-error-estimate} and \tw{Assumption}~\ref{the:elliptic-reconstruction-error-estimates}, respectively, it
  will suffice to estimate $\trecres$, or $\erecres$, in terms of
  $\trecerr$ or $\erecerr$.  

 \begin{Lem}[time reconstruction error bounds]
   \label{lem:t_rec2}
   For~$\rangefromto{n}1N$, we have the explicitly computable bounds
   \begin{equation}
     \dgenorm[0,t_n]{\trecerr}^2
     =
     \sum_{m=1}^n  \constref[\timestep m,\tdeg m]{eqn:time-reconstruction:constant}^2
     \Normreg{\jumpatof{m-1}{\full U}}_{\elldom}^2,
   \end{equation}
   and
   \begin{equation}
     \Norm{\trecerr}_{\leb{\infty}(0,t_n;\pivot)}^2= \max_{m=1,\dots,n} 	\Normreg{\jumpatof{m-1}{\full U}}_{\pivot}^2.
   \end{equation}
 \end{Lem}
 \begin{proof}
The proof follows immediately from Proposition \ref{pro:time-reconstruction-error-estimate}.
 \end{proof}
\begin{Lem}[space reconstruction error bounds]
  \label{lem:e_rec}
  Assume the availability of elliptic a posteriori error bounds as per Assumption \ref{the:elliptic-reconstruction-error-estimates}. Then, for~$\rangefromto{n}1N$, we have the explicitly computable bound
  \begin{equation}
    \dgenorm[0,t_n]{\erecerr}^2
    \leq
    \sum_{m=1}^n \int_{I_m} \cE_{\elldom,\xespace_n}[\full u, \feellop[n] \full u+\spaceproj\trecfeta{U}'-\trecfeta{U}']^2\d t.
  \end{equation}\qed
\end{Lem}
We can now prove the first main result of this work.
\oldtext{
\begin{The}[OLD RESULT $\leb{2}(I;\elldom)$-norm \aposteriori error bound]%
  Let Assumption \ref{the:elliptic-reconstruction-error-estimates}
  hold. Then, for $n=1,\dots,N$, we have the following \aposteriori
  error bound
  \begin{equation}%
    \begin{split}
      &\alpha^{-1}\ltwonorm{e(t_n\twchanges{^-})}^2+	\dgenorm[0,t_n]{e}^2\\
      \leq &\ 10c_{\alpha,\beta}\sum_{m=1}^n \int_{I_m}\!\! \cE_{\elldom,\xespace_n}[\full u, \feellop[n] \full u+\spaceproj\trecfeta{U}'-\trecfeta{U}']^2\d t
      +5\alpha^{-2}(\oscind[n])^2\\
      &+8c_{\alpha,\beta}	\sum_{m=1}^n  \constref[\timestep m,\tdeg m]{eqn:time-reconstruction:constant}^2
      \Normreg{\jumpatof{m-1}{\full U}}_{\elldom}^2
      +2{\alpha}^{-1}\!\!\max_{m=1,\dots,n} \!\!	\Normreg{\jumpatof{m-1}{\full U}}_{\pivot}^2 +2\alpha^{-1}(\eta^{\rm in})^2,
    \end{split}
  \end{equation}
  upon setting $c_{\alpha,\beta}:=1+(\nicefrac{\beta}{\alpha})^2$, with the constants~$\alpha,\beta$ from~\eqref{eq:ellop-prop}, and  $\oscind[n]:=\dgdnorm[0,t_n]{\fullproj[]  f-f }, $ and $\eta^{\rm in}:=\ltwonorm{u_0-\spaceproj[0]u_0}$ to signify the data oscillation and initial error indicators, respectively. 
\end{The}
\begin{proof}
  From
  \eqref{eqn:reconstruction-pdr:pointwise-form} and the PDE
  \eqref{eqn:def:exact-pde}, on each $I_n$, we have
	\begin{equation}%
	\pdt{\trecres\,}
	+\ellop\erecres
	=
	\pdt u+\ellop u
	-
	\pdt{\trecfeta u}
	\twchanges{-}
	\ellop\erec
	=
	f-\fullproj f.
	\end{equation}
	Testing \eqref{eqn:error_equation} with
	$\trecres$, we deduce
	\begin{equation}%
	\begin{split}
	\ltwop{ \pdt{\trecres\,}}{\trecres}
	+\duality{\ellop\erecres}{\erecres\, }
	= &\ 
	\duality{ f-\fullproj f}{\trecres\,}+\duality{\ellop\erecres}{\trecerr-\erecerr\, }
	.
	\end{split}
	\end{equation}
	Note that the term $\fullproj f -f$ occurring in
	the above estimator is computable
	signifying the so-called \emph{oscillation error}.
	Thanks to the continuity of $\trecres$, we can integrate the last identity on $(0,t_n)$, for some $n=1,\dots,N$, and through standard arguments, we obtain
	\begin{equation}%
	\begin{split}
	\frac12\ltwonorm{\trecres(t_n)}^2
	&+
	\alpha\dgenorm[0,t_n]{\erecres}^2\\
	\leq&\ 
	\frac12\ltwonorm{\trecres(0)}^2 +
	\beta
	\dgenorm[0,t_n]{\erecres}
	\dgenorm[0,t_n]{\trecerr-\erecerr}\\
	&+\dgdnorm[0,t_n]{f-\fullproj[] f}
	\big(\dgenorm[0,t_n]{\erecres}+\dgenorm[0,t_n]{\trecerr-\erecerr}\big),
	\end{split}
	\end{equation}
	using the continuity and coercivity of $\ellop$, cf.~\eqref{eq:ellop-prop}. Now, repeated use of the elementary inequality $2ab\le \epsilon a^2+\epsilon^{-1}b^2$ (for different values of $\epsilon>0$), and rearranging terms,
we deduce
	\begin{equation}%
	\begin{split}
\frac12\ltwonorm{\trecres(t_n)}^2+	\frac{\alpha}{2}\dgenorm[0,t_n]{\erecres}^2
	\leq&\ 
	\frac12\ltwonorm{\trecres(0)}^2 +
\frac{\alpha^2+\beta^2}{\alpha}
	\dgenorm[0,t_n]{\trecerr-\erecerr}^2\\
	&+\frac{5}{4\alpha}\dgdnorm[0,t_n]{f-\fullproj[] f}^2,
	\end{split}
	\end{equation}
	which implies
	\begin{equation}%
	\begin{split}
\frac{2}{\alpha}\ltwonorm{\trecres(t_n)}^2+	2\dgenorm[0,t_n]{\erecres}^2
\leq&\ 
4\Big(1+\frac{\beta^2}{\alpha^2}\Big)
\dgenorm[0,t_n]{\trecerr-\erecerr}^2\\
&+\frac{5}{\alpha^2}\dgdnorm[0,t_n]{f-\fullproj[] f}^2+	\frac2\alpha\ltwonorm{\trecres(0)}^2.
\end{split}
	\end{equation}
	Upon observing the estimates \begin{equation}
	\dgenorm[0,t_n]{e}^2\le 2\dgenorm[0,t_n]{\erecres}^2+2\dgenorm[0,t_n]{\erecerr}^2\end{equation}
	and
	\begin{equation}
\alpha^{-1}\ltwonorm{e(t_n\twchanges{^-})}^2\le \frac{2}{\alpha}\ltwonorm{\trecres(t_n)}^2+\frac{2}{\alpha}\ltwonorm{\trecerr(t_n\twchanges{^{-}})}^2,\end{equation}
 invoking Lemmas \ref{lem:t_rec2} and \ref{lem:e_rec},  and noting that $\ltwonorm{\trecres(0)} =\ltwonorm{u_0-\spaceproj[0]u_0}$, the result already follows.
\end{proof}
}

\newtext{
  \begin{The}[$\leb{2}(I;\elldom)$-norm \aposteriori error bound]
    \label{thm:L2H1}
    Let Assumption \ref{the:elliptic-reconstruction-error-estimates}
    hold. Then, for $n=1,\dots,N$,
    \tw{and 
      \begin{equation}\label{eq:theta}
	\left(1+\frac{\alpha}{\beta}\right)^{-1}\le\vartheta^2\le1+\frac{\alpha}{\beta},
      \end{equation}
    } 
    we have the following \aposteriori
    error bound
    \tw{
      \begin{equation}\label{eqn:L2_H1_apost}
	\begin{split}
	  \sqrt{\alpha}\dgenorm[0,t_n]{e}
	  \leq&\ 
	  \eta^{\rm in}
	  +
	  c_{\alpha,\beta}\vartheta^{-1}\oscind[n]
	  +c_{\alpha,\beta}\vartheta\beta
	  \left(\sum_{m=1}^n  \constref[\timestep m,\tdeg m]{eqn:time-reconstruction:constant}^2
          \Normreg{\jumpatof{m-1}{\full U}}_{\elldom}^2\right)^{\nicefrac12}\\&+\left(\sqrt{\alpha}+c_{\alpha,\beta}\vartheta\beta\right)
          \left(\sum_{m=1}^n \int_{I_m} \cE_{\elldom,\xespace_n}[\full u, \feellop[n] \full u+\spaceproj\trecfeta{U}'-\trecfeta{U}']^2\d t)\right)^{\nicefrac12},
	\end{split}
      \end{equation}
      where $c_{\alpha,\beta}:=\sqrt{\alpha^{-1}+\beta^{-1}}$, }
    with the constants~$\alpha,\beta$ from~\eqref{eq:ellop-prop}, and  $\oscind[n]:=\dgdnorm[0,t_n]{\fullproj[]  f-f }, $ and $\eta^{\rm in}:=\ltwonorm{u_0-\spaceproj[0]u_0}$ to signify the data oscillation and initial error indicators, respectively. 
  \end{The}
  \begin{proof}
    From
    \eqref{eqn:reconstruction-pdr:pointwise-form} and the PDE
    \eqref{eqn:def:exact-pde}, on each $I_n$, we have
    \begin{equation}\label{eqn:error_equation}
      \pdt{\trecres\,}
      +\ellop\erecres
      =
      \pdt u+\ellop u
      -
      \pdt{\trecfeta u}
      \twchanges{-}
      \ellop\erec
      =
      f-\fullproj f.
    \end{equation}
    Testing \eqref{eqn:error_equation} with
    $\trecres$, we deduce
    \begin{equation}\label{eqn:error_equation_new}
      \begin{split}
	\ltwop{ \pdt{\trecres\,}}{\trecres}
	+\duality{\ellop\erecres}{\erecres\, }
	= &\ 
	\duality{ f-\fullproj f}{\trecres\,}+\duality{\ellop\erecres}{\trecerr-\erecerr\, }
	.
      \end{split}
    \end{equation}
    Note that the term $\fullproj f -f$ occurring in
    the above estimator is computable
    signifying the so-called \emph{oscillation error}.
    Thanks to the continuity of $\trecres$, we can integrate the last identity on $(0,t_n)$, for some $n=1,\dots,N$, and through standard arguments, we obtain
    \begin{equation}\label{eqn:first_bound}
      \begin{split}
	\frac12\ltwonorm{\trecres(t_n)}^2
	+
	\alpha\dgenorm[0,t_n]{\erecres}^2
	&\leq\ 
	\frac12\ltwonorm{\trecres(0)}^2 +
	\beta
	\dgenorm[0,t_n]{\erecres}
	\dgenorm[0,t_n]{\trecerr-\erecerr}\\
	&\quad+\oscind[n]
	\big(\dgenorm[0,t_n]{\erecres}+\dgenorm[0,t_n]{\trecerr-\erecerr}\big),
      \end{split}
    \end{equation}
    using the continuity and coercivity of $\ellop$, cf.~\eqref{eq:ellop-prop}. 
    \tw{For any $\vartheta>0$
      observe the identity
      \begin{equation}
	\begin{split}
	  \beta
	  &\dgenorm[0,t_n]{\erecres}
	  \dgenorm[0,t_n]{\trecerr-\erecerr}
	  +\oscind[n]
	  \big(\dgenorm[0,t_n]{\erecres}+\dgenorm[0,t_n]{\trecerr-\erecerr}\big)\\
	  &=
	  \frac{\alpha}{2}\dgenorm[0,t_n]{\erecres}^2
	  +\frac12\left(\alpha^{-1}+\beta^{-1}\right)\left(\vartheta\beta\dgenorm[0,t_n]{\trecerr-\erecerr}+\vartheta^{-1}\oscind[n]\right)^2\\
	  &\quad-\frac{1}{2\alpha}\left(\beta\dgenorm[0,t_n]{\trecerr-\erecerr}
	  -\alpha\dgenorm[0,t_n]{\erecres}
	  +\oscind[n]\right)^2\\
	  &\quad+\frac12\left(\alpha^{-1}-{\vartheta^2}\left(\alpha^{-1}+\beta^{-1}\right)\right)\beta^2\dgenorm[0,t_n]{\trecerr-\erecerr}^2\\
	  &\quad+\frac12\left(\alpha^{-1}-{\vartheta^{-2}}\left(\alpha^{-1}+\beta^{-1}\right)\right)(\oscind[n])^2.
	\end{split}
      \end{equation}
      Notice that the last two terms are non-positive for $\vartheta$ as in~\eqref{eq:theta}.
      Then, from~\eqref{eqn:first_bound}, we infer
      \begin{equation}\label{eqn:L2_H1_one}
	\begin{split}
	  \alpha\dgenorm[0,t_n]{\erecres}^2
	  \le\ltwonorm{\trecres(0)}^2
	  +\left(\alpha^{-1}+\beta^{-1}\right)\left(\vartheta\beta\dgenorm[0,t_n]{\trecerr-\erecerr}+\vartheta^{-1}\oscind[n]\right)^2.
	\end{split}
      \end{equation}
      Taking square root and using the triangle inequality,
      $
      \dgenorm[0,t_n]{e}\le \dgenorm[0,t_n]{\erecres}+\dgenorm[0,t_n]{\erecerr}$, 
      we deduce
      \begin{equation}\label{eqn:L2_H1}
	\begin{split}
	  \sqrt{\alpha}\dgenorm[0,t_n]{e}
	  \leq&\ 
	  \ltwonorm{\trecres(0)} +c_{\alpha,\beta}\vartheta\beta
	  \dgenorm[0,t_n]{\trecerr}\\&+\left(\sqrt{\alpha}+c_{\alpha,\beta}\vartheta\beta\right)
	  \dgenorm[0,t_n]{\erecerr}+
	  c_{\alpha,\beta}\vartheta^{-1}\oscind[n].
	\end{split}
      \end{equation}
    }
    Invoking now Lemmas \ref{lem:t_rec2} and \ref{lem:e_rec},  and noting that $\ltwonorm{\trecres(0)} =\ltwonorm{u_0-\spaceproj[0]u_0}$, the result already follows.
  \end{proof}
}

\tw{
  \begin{Obs}
    For the particular case $\ellop=-\Delta$, $\elldom=\sobh1_0(\Omega)$ and $\pivot=\leb2(\Omega)$, cf.~Example~\ref{exa:concrete-elliptic-operators}, note that the constants in~\eqref{eq:ellop-prop} can be chosen to be $\alpha=\beta=1$. Especially, selecting $\vartheta^2=\nicefrac12$, we have $c_{\alpha,\beta}\vartheta=1$ in~\eqref{eqn:L2_H1_apost}.
  \end{Obs}
}

\begin{Obs}[mesh change error via elliptic reconstruction]
  \label{rem:mesh-change}
  The elliptic reconstruction
  \eqref{eqn:ell_rec_mod} features a
  mesh-change type term. This is an important departure from the (standard) elliptic
    reconstruction proposed in
    \cite{LakkisMakridakis:06:article:Elliptic}; cf. also \cite[Definition 6.1]{CGM}. For instance,
  if $\cE_{\elldom,\xespace_n}$ from~\eqref{eq:ellapost} is the standard
  residual energy norm \aposteriori error estimator, the element
  residual includes the expression
  $\spaceproj\trecfeta{U}'-\trecfeta{U}'$, which is
  effectively a mesh-change term. Indeed, from
  \eqref{eq:time-reconstruction-formula} we have on $I_n$:
  \begin{equation}
    \begin{split}
      \spaceproj\trecfeta{U}'-\trecfeta{U}'
      =&\
      \spaceproj
      \qp[big]{\full u'+\chi_n(\jumpatof{n-1}{\full u})}
      -
      \qp[big]{\full u'+\chi_n(\jumpatof{n-1}{\full u})}
      \\
      =&\
      \chi_n(\spaceproj\jumpatof{n-1}{\full u}-\jumpatof{n-1}{\full u})
      =\chi_n(\full U(t_{n-1}^-) - \spaceproj \full U(t_{n-1}^-)),
    \end{split}
  \end{equation}
  from the commutativity of the spatial projection $\spaceproj$ and
  operations on the time variable.  Let now
  $\varkappa_n :I_n\to \mathbb{R}$ be the
  polynomial of degree $r_n$ representing the time lifting $\chi_n$ from~\eqref{eq:tl2},
  i.e. such that
  \begin{equation}
    \label{eqn:def:time-lifting-K}
    \chi_n(\jumpatof{n-1}{\ta w})(t)
    =
    \varkappa_n (t)\jumpatof{n-1}{\ta w},
  \end{equation}
  for all $t\in I_n$ and $\ta W \in \cY$; we refer to \cite[Lemmas 6
    \& 7]{SW10} or~\cite[Remark~1]{HolmWihler:15} for an explicit formula for
  $\varkappa_n $. Then, we can
  conclude
  \begin{equation}
    \label{eqn:mesh-change_repr}
    \begin{split}
      \spaceproj\trecfeta{U}'(t)-\trecfeta{U}'(t)
      =&\
      \varkappa_n (t)
      \big( \full U(t_{n-1}^-) - \spaceproj \full U(t_{n-1}^-)\big)
      ,\qquad t\in I_n,
    \end{split}
  \end{equation}	
  i.e.~the arbitrary order \dG-analogue to the classical mesh-change
  indicator. We note that the representation
  \eqref{eqn:mesh-change_repr} is particularly relevant in
  implementation as it can be used to efficiently realize the space
  reconstruction error bounds in Lemma \ref{lem:e_rec}. Finally, observing that~$\varkappa_n$ takes its maximum at~$t_{n-1}$, which follows from~\cite[Lemma~6]{SW10}, and recalling~\eqref{eq:tl2}, it is possible to compute the maximum value of
  $\varkappa_n$ on $\overline I_n$:
  \begin{equation}
  \Norm{\varkappa_n}^2_{\leb\infty(I_n)}
  =|\varkappa_n(t_{n-1})|^2
  =\int_{I_n}\varkappa_n(t)^2\d t.
  \end{equation}
  Hence
  \begin{equation}
  \Norm{\varkappa_n}_{\leb\infty(I_n)}
  =\Norm{\varkappa_n}_{\leb2(I_n)}
  =\frac{r_n+1}{\sqrt{\tau_n}},
  \end{equation}
  see~\cite[Proposition~2]{SW10}.
\end{Obs}

\begin{Obs}[reconstruction vs direct approach]
    \label{obs:reconstruction-vs-direct}
    As noted by one referee, it is possible to spare the use of
    elliptic reconstruction in the proof of the
    $\leb{2}(I;\elldom)$-norm \aposteriori bound, for particular cases
    of operators $\ellop$ and of Gelfand triples~\eqref{eq:Gelfand},
    as discussed in \cite{ESV:16,ESV:17, GKSZ17} for the \tw{specific}
    case $\ellop=-\Delta$, $\elldom=\sobh1_0(\Omega)$ and
    $\pivot=\leb2(\Omega)$;
    \tw{cf.~Example~\ref{exa:concrete-elliptic-operators}}. As a
    result, different \aposteriori error estimators arise with,
    possibly, slightly different constants multiplying common terms in
    the estimator. Elliptic reconstruction, nonetheless, offers the
    ability to use various elliptic estimators from the literature in
    the bound: this feature may become important for multiscale
    operators $\ellop$, e.g., singular perturbations. In general,
    elliptic reconstruction allows for, crucial in some cases,
    flexibility in the handling for more complicated spatial operators
    $\ellop$, e.g., nonlinear or singularly perturbed operators
    \cite{CGM, GM,mohammad}.
\end{Obs}
\section{$\leb \infty(I;\pivot)$-norm \aposteriori error analysis}
\label{sec:Linfty(L2)apost}
We continue by proving an \aposteriori error bound in the $\leb \infty(I;\pivot)$-norm, which appears to be of higher order than the $\leb{2}(I;\elldom)$-norm one presented above. A key difference with respect to the above proof of the $\leb 2(I;\elldom)$-norm bound is the use of a combined space-time reconstruction defined below.
\subsection{Error-residual relation}
\label{sec:error-residual-relation}
Set $w:=\widehat{\erec}$, i.e. the time reconstruction of the elliptic reconstruction, noting that Definition \ref{def:time-reconstruction} is still valid for functions on $\elldom$ in space. 
Now, introducing the \emph{time error} 
\begin{equation}\label{eq:rho1}
  \rho:=w-u,
\end{equation}
which is continuous on~$\closure I=[0,T]$, and subtracting the PDE \eqref{eqn:def:exact-pde} from~\eqref{eqn:reconstruction-pdr:pointwise-form}, we have: 
\begin{equation}
  (\trecfeta u-u)'
  +
  \ellop (\erec-u)=
  \fullproj[]  f-f,
\end{equation}
our aim being to deduce an evolution equation for $\rho$. To this end, we have
\begin{equation}\label{eq:spacetime_strong_two}
  \begin{split}
    \rho'
   +
   \ellop \rho=&\ 
   \fullproj[]  f -f+ (w-\trecfeta u)'+ \ellop (w-\erec)\\
   =&\ \fullproj[]  f-f + ( \widehat{\erec-\full u})'+ \ellop (w-\erec)=:\xi,
  \end{split}
\end{equation}
from the linearity of the time reconstruction. 

\newtext{
We define $\mini{\mesh[n]t}  {\mesh[n-1]t}$ to be the finest common coarsening of the two meshes, \tw{with} the associated largest common finite element subspace \tw{given by} $\fesm n=\fes n\meet\fes{n-1}$, and a meshsize $\meshsizem=\maxi{\meshsize[n]}{\meshsize[n-1]}$.
\begin{Lem}[time reconstruction error bound]
  \label{proposition_sigma_direct}
  Let $ t\in I_n$, $n=0,1\dots, N$. Then, we have the abstract estimate
  \begin{equation}\label{first_sigma_estimate}
    \dgdnorm[I_n]{
      \ellop(w-\erec)} \le \constref[\timestep n,\tdeg n]{eqn:time-reconstruction:constant} \eta_{n,1}^{\rm time},
  \end{equation}
  with
  \begin{equation}
  \eta_{n,1}^{\rm time}:=\| \jumpatof{n-1}{
    \fullproj[]  f-\pdt{\trecfeta u}}\|_{\dual\elldom}.
  \end{equation}
  Assume further that, for every $v\in\elldom$, there exists a $\fe v\in \xespace_n^\ominus$, such that
  \begin{equation}\label{assumptions_on_spaces}
    \|(\meshsizem)^{-s}(v-\fe v)\|_{\pivot}\le C_{\rm ap}\|v\|_{\elldom}, \quad\text{and} \quad \|\fe v\|_{\elldom}\le C_{\rm stab}\|v\|_{\elldom},
  \end{equation}
  for some $s>0$, and for generic constants $C_{\rm ap}, C_{\rm stab}>0$, independent of $v$ and of $\meshsizem$. Then, we also have the bound
  \begin{equation}\label{second_sigma_estimate}
    \dgdnorm[I_n]{
      \ellop(w-\erec)} \le \constref[\timestep n,\tdeg n]{eqn:time-reconstruction:constant} \eta_{n,2}^{\rm time},
  \end{equation}
  with
  \begin{equation}
  \eta_{n,2}^{\rm time}:= C_{\rm ap}\|(\meshsizem)^{s} \jumpatof{n-1}{
    \fullproj[]  f-\pdt{\trecfeta u}}\|_\pivot +C_{\rm stab}	\|\jumpatof{n-1}{\full u}\|_\elldom.
  \end{equation}
  Therefore, setting $\eta_n^{\rm time}:=\min\{\eta_{n,1}^{\rm time},\eta_{n,2}^{\rm time} \}$, we have the combined estimate 
  \begin{equation}
  \dgdnorm[I_n]{
    \ellop(w-\erec)} \le \constref[\timestep n,\tdeg n]{eqn:time-reconstruction:constant} \eta_{n}^{\rm time}.
  \end{equation}
\end{Lem}
}

\newtext{
\begin{proof} 	Recalling~\eqref{eq:time-reconstruction-formula} and using the fact
	that the elliptic operator~$\ellop$ is time independent, we
	immediately observe that 
	 \begin{equation}\label{eq:acommute}
	\ellop w=\widehat{\ellop\erec}.
	\end{equation}
	Therefore, by Proposition~\ref{pro:time-reconstruction-error-estimate} and of Remark \ref{rem:magic_identity}, for~$\rangefromto{n}1{N-1}$, we conclude that
	\begin{equation}\label{eq:413a}
	\begin{aligned}
	\constref[\timestep n,\tdeg n]{eqn:time-reconstruction:constant}^{-1}		\dgdnorm[I_n]{
		\ellop(w-\erec)}
	&=
	\Norm{\jumpatof{n-1}{\ellop\erec}}_{\dual\elldom}
	=
	\Norm{\jumpatof{n-1}{\fullproj[]f-\trecfeta u'}}_{\dual\elldom},
	\end{aligned}
	\end{equation}
	\tw{which is~\eqref{first_sigma_estimate}. Furthermore, 
exploit} the orthogonality identity
		\begin{equation}\label{orthog}
	\duality{\jumpatof{n-1}{\ellop\erec}}{\fe v}	= \ltwop{\jumpatof{n-1}{\feellop[n] \full u}}{\fe v}=	\duality{\jumpatof{n-1}{\ellop\full u}}{\fe v}\quad \text{for all } \fe v\in \xespace_n^{\ominus},
		\end{equation}
		which follows from the definition of $\erec$ and the properties of $\fullproj[]$. From the latter, together with Remark \ref{rem:magic_identity}, and the temporal independence of $\ellop$, we have
			\begin{equation}
		\begin{aligned}
		\duality{\jumpatof{n-1}{\ellop\erec}}{v}
		=&\  (\jumpatof{n-1}{
			\fullproj[]  f-\pdt{\trecfeta u}},v-\fe v)_\pivot +	\duality{\ellop\jumpatof{n-1}{\full u}}{\fe v}
		\quad \text{for all } \fe v\in \xespace_n^{\ominus}.
		\end{aligned}
		\end{equation}
		From \eqref{assumptions_on_spaces} and the continuity of $\ellop$, therefore, we deduce
		\begin{equation}
		\begin{aligned}
		\duality{\jumpatof{n-1}{\ellop\erec}}{v}
		\le &\  \big(C_{\rm ap}\|(\meshsizem)^{s} \jumpatof{n-1}{
			\fullproj[]  f-\pdt{\trecfeta u}}\|_\pivot +C_{\rm stab}	\|\jumpatof{n-1}{\full u}\|_\elldom\big) \| v\|_\elldom,
	\end{aligned}
\end{equation}
and the result already follows.
\end{proof}
}

\begin{Obs}[Practical upper bound for $\eta_{n,1}^{\rm time}$]
\newtext{In standard settings, e.g., in the canonical example of \tw{the} heat equation, \tw{cf.~Example~\ref{exa:concrete-elliptic-operators},} for which \eqref{assumptions_on_spaces} holds, we can use $\eta_{n,2}^{\rm time}$ for $\eta_{n}^{\rm time}$; this is particularly pertinent when the dimension of $\xespace_n^{\ominus}$ is comparable to those of $\xespace_n$ and $\xespace_{n-1}$. In cases, however, in which \eqref{assumptions_on_spaces} may not be effective, or even known, we may be forced to revert to \eqref{first_sigma_estimate}; see \cite{oliver} for two such settings, \tw{where} one \tw{involves} virtual element methods, and the other \tw{incorporates} moving mesh methods. The dual norm in $\eta_{n,1}^{\rm time}$ typically requires a global inversion to be evaluated, in the absence of a natural Galerkin orthogonality property like \eqref{orthog} above. 
}	
	
\newtext{
	Although solving an elliptic problem to determine/estimate the time estimator $\eta_{n,1}^{\rm time}$ may appear to be computationally demanding, such global \tw{solves} may be typically performed in a lower dimensional space than $\xespace_n$ itself, through an \aposteriori error controlled fashion, as we now demonstrate. Let $\varphi \in \elldom$ be the solution to the problem $\ellop \varphi = w$ \tw{(with $w:=\jumpatof{n-1}{\fullproj[]  f-\pdt{\trecfeta u}}$ for instance)}.
	Let also \tw{an} approximation $\Phi\in \tilde{\xespace}_n$ to $\varphi$, for some Galerkin space $\tilde{\xespace}_n$, \tw{be} given by 
\begin{equation}\label{ell_neg_norm}
	\duality{\ellop \Phi}{\fe v}=	\duality{w}{\fe v}\quad\text{ for all } \fe v\in \tilde{\xespace}_n.
\end{equation}	 
Now, the continuity of $\ellop$ implies
$
\|w\|_{\dual\elldom}\le \beta \|\varphi\|_{\elldom}.
$
	\tw{In addition}, from coercivity, continuity and the Galerkin orthogonality $\duality{\ellop (\varphi- \Phi)}{\Phi}=0 $, we have, respectively,
	\begin{equation}
	\alpha \|\varphi\|_{\elldom}^2 \le \duality{\ellop \varphi}{\varphi}= \duality{\ellop (\varphi- \Phi)}{\varphi- \Phi} +\duality{\ellop \Phi}{\Phi}\le 
	\beta\big( \|\varphi-\Phi\|_{\elldom}^2+ \|\Phi\|_{\elldom}^2\big),
	\end{equation}
	using also the self-adjointness of $\ellop$. The above, together with the assumed availability of an \aposteriori error indicator in the energy norm give  the computable bound
	\begin{equation}\label{computable_full_solve}
\|w\|_{\dual\elldom}^2\le \frac{\beta^3}{\alpha} 
	\big(\cE_{\elldom,\tilde{\xespace}_n}^2[\Phi,w]+ \|\Phi\|_{\elldom}^2\big),
	\end{equation}
	or, for the specific case $ w=\jumpatof{n-1}{
		\fullproj[]  f-\pdt{\trecfeta u}} $,
	\begin{equation}
	\eta_{n,1}^{\rm time}\le \sqrt{\nicefrac{\beta^3}{\alpha}} \big(\cE_{\elldom,\tilde{\xespace}_n}^2[\Phi,\jumpatof{n-1}{
		\fullproj[]  f-\pdt{\trecfeta u}}]+ \|\Phi\|_{\elldom}^2\big)^{\nicefrac12}.
	\end{equation}
}

\oldtext{
REMOVE---For~$\rangefromto{n}{1}{N}$, we introduce the
  oscillation (in time) of the spatial $\pivot$-projection $\spaceproj f$ of $f$ about its space-time $\pivot$-projection onto the time-discrete Galerkin space:
\begin{equation}
    \timeosc f|_{I_n}:=\fullproj f-\spaceproj f.
  \end{equation}
Owing to the triangle inequality, it holds that
\begin{align}\label{eq:timeindtriangle}
\eta_{n,1}^{\rm time}
\le\Norm{\jumpatof{n-1}{\timeosc f-\trecfeta u'}}_{\dual\elldom}+\Norm{\jumpatof{n-1}{\fullproj[]f-\timeosc f}}_{\dual\elldom}.
\end{align}
Observe that both $\jumpatof{n-1}{\timeosc f}$ and
$\jumpatof{n-1}{\trecfeta u'}$ are temporal oscillations/approximations. Therefore, it is possible that, bounding the first term on the right-hand side of
\eqref{eq:timeindtriangle} by a Poincar\'e--Friedrichs inequality, viz.,
\begin{equation}
  \label{eqn:poincare-based-jump-estimate}
  \Norm{\jumpatof{n-1}{\timeosc f-\trecfeta u'}}_{\dual\elldom}
  \le
  \sqrt{\constext{PF}}
  \Norm{\jumpatof{n-1}{\timeosc f-\trecfeta u'}}_{\pivot},
\end{equation}
for some positive constant $\constext{PF}$ independent of the functions involved, may not be detrimental to the convergence rate. For the second term on the right-hand side of \eqref{eq:timeindtriangle}, provided that the one-sided limits~$f(t_{n-1}^\pm)$ exist in~$\dual\elldom$, we notice that
     \begin{equation}
        \label{eq:concrete-Linfty(I;pivot)-timeindicator:splitting}
        \begin{split}
          \jumpatof{n-1}{\fullproj[]f-\timeosc f}
          &
          =
          \spaceproj f (t_{n-1}^+)-\spaceproj[n-1] f(t_{n-1}^-)
          \\
          &
          =
          \spaceproj f (\tp{n-1})
          -
          \spaceproj[n] f (\tm{n-1})
          +
          \spaceproj[n] f (\tm{n-1})
          -
          \spaceproj[n-1] f(\tm{n-1})
          \\
          &=
          \memotempa
          +
          \memotempb
          ;
        \end{split}
      \end{equation}%
      here $\spaceproj[ ]$ indicates the $\cI$-piecewise constant time-extension
      of $\spaceproj$, for $n\integerbetween1N$.
In particular, if the Galerkin space does \emph{not} change across~$t_{n-1}$, i.e. if it holds~$\fes{n-1}=\fes{n}$, and if  $\jumpatof{n-1}{f}=0$ in~$\dual\elldom$, 
then $\jumpatof{n-1}{\fullproj[]f-\timeosc f}=0$, and thus, \eqref{eq:timeindtriangle} simplifies to
      \begin{equation}
        \begin{split}
          \timeind
          &\le
          \constext{PF}^{\nicefrac12}
          \Norm{\jumpatof{n-1}{\timeosc f-\trecfeta u'}}_{\pivot}
          ,
        \end{split}
      \end{equation}
which is easily computable.}
\end{Obs}

\oldtext{
REMOVE---
\begin{Obs}[computable upper bound for $\timeind$ in the context of spatial FEM]
We address now the practical implications of the computation of~\eqref{eq:concrete-Linfty(I;pivot)-timeindicator:splitting}
    when $\pivot=\leb2(\W)$ and $f(t)\in\pivot$ for $0\leq t\leq T$, to focus on most concrete situations where this is the case.
      More specifically, for~$\rangefromto{n}{1}{N}$, consider a classical finite element discretization space 
      \begin{equation}
      \fes n=\poly
      p(\mesh[n]t)\meet\elldom,
      \end{equation}
      for a conforming and shape-regular partition $\mesh[n]t$ 
      of meshsize (function) $\meshsize[n]$ of the underlying spatial domain, and a (fixed) polynomial degree~$p\ge 1$.
      First, we note that the term $\memotempa$ appearing in~\eqref{eq:concrete-Linfty(I;pivot)-timeindicator:splitting} is
      present only if the function $f$ is discontinuous
      across $\t{n-1}$. In that case it may be bounded in various ways depending on the respective application context. For instance,
      since $\memotempa\in\xespace_n$, it is computable
      and belongs to $\leb2{(\W)}$ and, thus, we have
      \begin{equation}
        \duality{\memotempa}{v}
        =
        (\memotempa,v)_{\leb2{(\W)}}
        \leq
        \constext{PF}^{\nicefrac12}
        \Norm{\memotempa}_{\leb2{(\W)}}
        \Norm v_{\elldom},
      \end{equation}
      owing to~\eqref{eq:identify} and the Poincar\'e--Friedrichs inequality (which would not degrade
      the bound because $\memotempa$ is only related to the temporal smoothness
      of $f$). Alternatively, if it is deemed appropriate the $\dual\elldom$-norm of $\memotempa$ can be computed by means of an additional elliptic solve involving
      $\feellop$.
	  To deal with~$\memotempb$, we define $\mini{\mesh[n]t}  {\mesh[n-1]t}$ to be the finest common refinement of the two meshes, and consider the associated largest common finite element subspace $\fesm n=\fes n\meet\fes{n-1}$ of $\fes{n-1}$ and $\fes n$, and a meshsize $\meshsizem=\maxi{\meshsize[n]}{\meshsize[n-1]}$.    
      Next, denoting by $\clementprojm$ a Cl\'ement-type
      interpolation operator onto $\fesm n$, for each $v\in\elldom$ we have,
      \begin{equation}
        \begin{split}
          \duality{
            \memotempb
          }v
          &
          =
          (\memotempb,v)_{\leb2{(\W)}}
          \\
          &
          =
          (\memotempb,v-\clementprojm v)_{\leb2{(\W)}}
          \\
          &
          \leq
          \constext{Cl}
          \Norm{
            \meshsizem
            \memotempb
          }_{\leb2(\W)}
          \Norm{v}_\elldom,
        \end{split}
      \end{equation}
      with an interpolation constant $\constext{Cl}>0$.
Then, we infer that
      \begin{equation}
        \begin{split}
        \Norm{\memotempb}_{\dual\elldom}
          \leq
          \constext{Cl}
          \Norm{
            \meshsizem
            \memotempb
          }_{\leb2(\W)}
          ,
        \end{split}
      \end{equation}
      which is a computable estimate. 
\end{Obs}
}

\begin{Lem}[space reconstruction error bound]
  \label{lem:ell_rec}
  For $t\in \timeinter n$, $\rangefromto n1N$, we have the bound
  \begin{equation}
    \label{eq:lem:ell_rec-inequality}
    \Norm{(\widehat{\erec-\full u})'(t)}_{\dual\elldom%
    }
    \le
    \spaceind
    ,
  \end{equation}
  where
  \begin{equation}
    \begin{split}
      \spaceind (t):=&\
      \cE_{\dual\elldom,\xespace_n}[{
          \full u'(t), (\feellop \full u' +\spaceproj\trecfeta{U}''-\trecfeta{U}'')(t)
      }]
      \\
      &
      +
      |\varkappa_n (t)|
      \cE_{\dual\elldom,\xespace_n}[\full u(t_{n-1}^+),	(\feellop 
	\full u+ \spaceproj\trecfeta{U}'-\trecfeta{U}')(t_{n-1}^+)]
      \\
      &+
      |\varkappa_n (t)|\cE_{\dual\elldom,\xespace_{n-1}}[\full u(t_{n-1}^-),
        (\feellop [n-1] \full u
        +
        \spaceproj[n-1]\trecfeta{U}'
        -
        \trecfeta{U}')(t_{n-1}^-)],
    \end{split}
  \end{equation}
  with $\varkappa_n :I_n\to\mathbb{R}$ from
  \eqref{eqn:def:time-lifting-K} in Remark \ref{rem:mesh-change}, and
  the estimator is explicitly computable, via Assumption
  \ref{the:elliptic-reconstruction-error-estimates}.
\end{Lem}

\begin{proof} Recalling~\eqref{eq:time-reconstruction-formula}, on~$I_n$ we have 
  \begin{equation}
    ( \widehat{\erec-\full u})'=
    \pdt\erec+\chi_n\qpreg{\jumpatof{n-1}{\erec}}-\pdt{\full u	}
    -\chi_n\qp{\jumpatof{n-1}{\full u}}.
  \end{equation}
  Hence, with the notation of Remark \ref{rem:mesh-change}, it follows that
  \begin{equation}\label{eq:aux20190603}
  	\Norm{(\widehat{\erec-\full u})'(t)}_{\dual\elldom}
	\le
	\Norm{\erec'-\full u'}_{\dual\elldom}
	+
	|\varkappa_n(t)|\Norm{\jumpatof{n-1}{\erec-\full u}}_{\dual\elldom},\quad t\in I_n.
  \end{equation}
  Differentiating \eqref{eqn:erec_orthog} with respect to $t\in I_n$,
  and recalling that the elliptic operator~$\ellop$ is $t$-independent, we deduce the Galerkin orthogonality relation
  \begin{equation}
    \duality{\ellop
      \pdt\erec-\ellop
      \pdt{\full u}
    }{\fe v}
    = 0
    \Foreach
    \fe v\in\fes n.
  \end{equation}
 We conclude that $\pdt{\full u} \in\fes n$ is the Galerkin approximation of~$\erec'$, where
  \begin{equation}\label{eqn:ell_rec_one}
    \ellop\erec'
    =
    \feellop 
    \full u'+ \spaceproj\trecfeta{U}''-\trecfeta{U}''
    \quad  \text{ for }
    t\in I_n,
  \end{equation}
   by differentiating
  \eqref{eqn:def:elliptic-reconstruction:strong} with respect to $t$.
  Furthermore, from the definition of the elliptic reconstruction~\eqref{eqn:def:elliptic-reconstruction:strong}, we have, respectively,
  \begin{equation}\label{eqn:gal_proj_two}
    \begin{split}
      \ellop\erec(t_{n-1}^+)
      =&\ 
      \feellop 
      \full u(t_{n-1}^+)+ \spaceproj\trecfeta{U}'(t_{n-1}^+)-\trecfeta{U}'(t_{n-1}^+),\\
      \ellop\erec(t_{n-1}^-)
      =&\ 
      \feellop [n-1]
      \full u(t_{n-1}^-)
      +
      \spaceproj[n-1]\trecfeta{U}'(t_{n-1}^-)-\trecfeta{U}'(t_{n-1}^-),
    \end{split}
  \end{equation}
  implying, therefore,
  \begin{equation}
    \begin{split}
      \duality{	\ellop\erec(t_{n-1}^+)-\ellop\full u(t_{n-1}^+)}{\fe v}  =&\  0
      \quad  \Foreach
      \fe v\in\fes n,
      \\
      \duality{	\ellop\erec(t_{n-1}^-)-\ellop\full u(t_{n-1}^-)}{\fe v}  =&\  0
      \quad  \Foreach
      \fe v\in\fes {n-1}.
    \end{split}
  \end{equation}
  This means that $\full u(t_{n-1}^+)\in \fes n$ and $\full
  u(t_{n-1}^-)\in \fes {n-1}$ are the Galerkin approximations of the first and second problem in \eqref{eqn:gal_proj_two}, respectively.
  
Now, \eqref{eqn:ell_rec_one} together with Assumption
  \ref{the:elliptic-reconstruction-error-estimates}, yield
  \begin{equation}
    \Norm{\erec'  -\full u'}_{\dual\elldom}
    \leq
    \cE_{\dual\elldom,\xespace_n}[{
        \full u', \feellop \full u' +\spaceproj\trecfeta{U}''-\trecfeta{U}''
    }]\quad\text{on }I_n.
  \end{equation}
  Similarly, noting that
  \begin{equation}\label{eqn:ell_gal_jump}
    \begin{split}
      \jumpatof{n-1}{\erec-\full u}
      =
        (\erec-\full u)(t_{n-1}^+)-(\erec-\full u)(t_{n-1}^-)
      ,
    \end{split}
  \end{equation} 
  and combining \eqref{eqn:gal_proj_two} with Assumption
  \ref{the:elliptic-reconstruction-error-estimates}, leads to
  \begin{equation}
    \begin{split}
      \Norm{	\jumpatof{n-1}{\erec-\full u}}_{\dual\elldom}
      \leq&\
      \cE_{\dual\elldom,\xespace_n}[\full u(t_{n-1}^+),	(\feellop 
        \full u+ \spaceproj\trecfeta{U}'-\trecfeta{U}')(t_{n-1}^+)]\\
      &+
      \cE_{\dual\elldom,\xespace_{n-1}}[{
          \full u(t_{n-1}^-),
          ({
            \feellop [n-1] \full u
            + \spaceproj[n-1]\trecfeta{U}'-\trecfeta{U}'
          })(t_{n-1}^-)
      }].
    \end{split}
  \end{equation}
  The result already follows by inserting the above estimates into~\eqref{eq:aux20190603}.
\end{proof}

We are now ready to present the second main result of this work.
\begin{The}[$\leb{\infty}(I;\pivot)$-norm \aposteriori bound]
  \label{the:LinfL2final}
  With the notation of \S\ref{sec:setup} and under
  Assumption~\ref{the:elliptic-reconstruction-error-estimates}, along
  with the assumptions of Lemmata \ref{proposition_sigma_direct} and
  \ref{lem:ell_rec}, for each $\rangefromto n1N$, we
  have the bound
  \begin{equation}\label{eq:Linfty_final}
    \begin{aligned}
      \Norm{u-\full u}&_{\leb{\infty}(0,t_n;\pivot)}\\
      &\le
     \eta^{\rm in}
      +
      (2\alpha)^{-\nicefrac12} \oscind[n]+\max_{j=1,\ldots,n} 
      \Norm{\jumpatof{j-1}{\full u}}_{\pivot}
      \\
      &
      \quad
      +
      \qp[Big]{
        \frac{1}{2\alpha}\sum_{m=1}^n
        \constref[\timestep m,\tdeg m]{eqn:time-reconstruction:constant}^2
        (\timeind[m])^2
      }^{\nicefrac12} %
      + \qp[Big]{\frac{1}{2\alpha}
      \sum_{m=1}^n \int_{I_m}(\spaceind[m])^2\d t
      }^{\nicefrac12}
      \\
      &\quad
      +3\max_{j=1,\ldots,n}
      \sup_{%
        I_j}
      \cE_{\pivot,\fes j}[{
          \full u, \feellop[j] \full u +\spaceproj[j]\trecfeta{U}'-\trecfeta{U}'
      }],
    \end{aligned}
  \end{equation}    
  with $\oscind[n]$ and $\eta^{\rm in}$ as in Theorem \ref{thm:L2H1}.
\end{The}
\begin{proof}
  We start from \eqref{eq:spacetime_strong_two}, which upon testing
  with $\rho$ and integrating with respect to $t$ between $(0,t)$,
  along with the coercivity of $\ellop$, \eqref{eq:ellop-prop}, gives
  \begin{equation}	\label{eq:Linfty_mid0}
    \frac{1}{2}\ltwonorm{{\rho}(t)}^2 +\alpha \dgenorm[0,t]{\rho}^2
    \le
    \frac{1}{2}\ltwonorm{{\rho}(0)}^2+\int_0^t \duality{\xi}{\rho}\d s,
  \end{equation}    
  and, thus,
  \begin{equation}
    \label{eq:Linfty_mid1}
    \ltwonorm{{\rho}(t)}^2 
    \le
    \ltwonorm{{\rho}(0)}^2+\frac{1}{2\alpha} \dgdnorm[0,t]{\xi}^2,
  \end{equation}    
  via standard arguments.
  
  Now fix $n\in\{1,\dots, N\}$ and choose $t^\star\in [0,t_n]$ such that 
  \begin{equation*}
    \ltwonorm{{\rho}(t^\star)}
    =
    \Norm{\rho}_{\leb{\infty}(0,t_n;\pivot)}.\end{equation*}
  Then, letting $t=t^\star$ in~\eqref{eq:Linfty_mid1}, and taking the square root, yields 
  \begin{equation}
    \label{eq:Linfty_mid2}
    \Norm{\rho}_{\leb{\infty}(0,t_n;\pivot)}
    \le
    \ltwonorm{{\rho}(0)}+(2\alpha)^{-\nicefrac12} \dgdnorm[0,t_n]{\xi}.
  \end{equation}    
  It remains to bound the last term on the right-hand side of
  \eqref{eq:Linfty_mid2}. To that end, applying~\eqref{eq:rec1st}, and involving Lemmata~\ref{proposition_sigma_direct}
  and~\ref{lem:ell_rec}, results in
  \begin{equation}
    \label{eqn:bound_on_xi}
    \begin{split}
      \dgdnorm[0,t_n]{\xi}\le&\  
      \oscind[n]+\dgdnorm[0,t_n]{\ellop( w-\erec)}
      +
      \dgdnorm[0,t_n]{( \widehat{\erec-\full u})'}
      \\
      \le &\ \oscind[n]+\!
      \Big(
      \sum_{m=1}^n\constref[\timestep m,\tdeg m]{eqn:time-reconstruction:constant}^2
      (\timeind[m])^2\Big)^{\nicefrac12} \!\!
      +\!
      \Big(
      \sum_{m=1}^n \int_{I_m}(\spaceind[m])^2\d t
      \Big)^{\nicefrac12}\!.
    \end{split}
  \end{equation}
  
  The triangle inequality now gives
  \begin{equation}
  \begin{split}\label{eq:tri2}
    \Norm{u-\full u}&_{\leb{\infty}(0,t_n;\pivot)}\\
    &\le 
    \Norm{{\rho}}_{\leb{\infty}(0,t_n;\pivot)}
    +\Norm{w - \erec}_{\leb{\infty}(0,t_n;\pivot)}
    +\Norm{\erec-\full u}_{\leb{\infty}(0,t_n;\pivot)}.
    \end{split}
  \end{equation}
  To estimate the second and third term on the right-hand side
  of~\eqref{eq:tri2}, we apply~\eqref{eq:time-reconstruction-linf}, and the
  triangle inequality, to obtain
    \begin{equation}
      \begin{aligned}
	\Norm{w - \erec}_{\leb{\infty}(0,t_n;\pivot)}
	&=
	\max _{j=1,\dots,n} \Norm{\jumpatof{j-1}{\erec}}_{\pivot}
	\\
	&\le
	\max _{j=1,\dots,n} \left(\Norm{\jumpatof{j-1}{\erec-\full u}}_{\pivot}
	+
	\Norm{\jumpatof{j-1}{\full u}}_{\pivot}\right).
      \end{aligned}
    \end{equation}
    Also, with the aid of~\eqref{eq:elliptic_apost},
    \begin{equation}
      \Norm{\erec-\full u}_{\leb{\infty}(0,t_n;\pivot)}
      \le{
        \max_{j=1,\ldots,n}\sup_{t\in I_j} 
	\cE_{\pivot,\fes j}[
          \full u, \feellop[j] \full u +\spaceproj[j]\trecfeta{U}'-\trecfeta{U}'
        ]
      }
      .
    \end{equation}
    Combining the last two estimates, we conclude
    \begin{equation}
      \begin{split}\label{eqn:ell_t_error_linftyl2} 
	&\Norm{w - \erec}_{\leb{\infty}(0,t_n;\pivot)}
	+\Norm{\erec-\full u}_{\leb{\infty}(0,t_n;\pivot)}\\
	\le &\  	{
          \max_{j=1,\ldots,n} \Big(3\sup_{t\in I_j}
	  \cE_{\pivot,\fes j}[
            \full u, \feellop[j] \full u +\spaceproj[j]\trecfeta{U}'-\trecfeta{U}'
        ]} 
	+
	\Norm{\jumpatof{j-1}{\full u}}_{\pivot}\Big).
      \end{split}
    \end{equation}
    Combining the above completes the argument.
\end{proof}
\begin{Obs}[alternative `short-time' \aposteriori error estimator]
  \label{obs:alternative-short-time-estimator}
  Starting from the alternative estimate
  \begin{equation}
    \int_0^t \duality{\xi}{\rho}\d s
    \le
    \Norm{\xi}_{\leb{1}(0,t_n;\pivot)}\Norm{\rho}_{\leb{\infty}(0,t_n;\pivot)},
  \end{equation}
  the bound~\eqref{eq:Linfty_mid0} with $t=t^\star$ as in the proof of Theorem \ref{the:LinfL2final} gives
  \begin{equation}
    \Norm{\rho}_{\leb{\infty}(0,t_n;\pivot)}^2
    \le
    2\ltwonorm{{\rho}(0)}^2+4 \Norm{\xi}_{\leb{1}(0,t_n;\pivot)}^2.
  \end{equation}    
  Triangle inequality trivially yields
  \begin{equation}
    \begin{split}
      \Norm{\xi}_{\leb{1}(0,t_n;\pivot)}\le&\  
      \Norm{\fullproj[]  f-f }_{\leb{1}(0,t_n;\pivot)}
      +
      \Norm{\ellop( w-\erec)}_{\leb{1}(0,t_n;\pivot)}\\
      &+
      \Norm{( \widehat{\erec-\full u})'}_{\leb{1}(0,t_n;\pivot)}.
    \end{split}
  \end{equation}
  Now, using H\"older's inequality, recalling~\eqref{eq:acommute}, and applying Proposition~\ref{pro:time-reconstruction-error-estimate} and
  Remark~\ref{rem:magic_identity}, we have, respectively,
  \begin{equation}\label{eq:413}
    \begin{aligned}
      \Norm{
        \ellop(w-\erec)}_{\leb{1}(0,t_n;\pivot)}^2	
      & \le t_n\Norm{	\ellop(w-\erec)}_{\leb{2}(0,t_n;\pivot)}^2
      \\
      &=t_n\sum_{m=1}^n
      \constref[\timestep m,\tdeg m]{eqn:time-reconstruction:constant}^2
      \Norm{\jumpatof{m-1}{\ellop\erec}}_{\pivot}^2\\
      &=t_n
      \sum_{m=1}^n\constref[\timestep m,\tdeg m]{eqn:time-reconstruction:constant}^2
      \Norm{\jumpatof{m-1}{\fullproj[]f-\trecfeta u'}}_{\pivot}^2.
    \end{aligned}
  \end{equation}
  In addition, an inspection of the proof of Lemma \ref{lem:ell_rec} reveals the bound
  \begin{equation}
    \Norm{( \widehat{\erec-\full u})'(t)}_{\pivot}
    \le
    \spaceind[m]^\star(t)
    ,
  \end{equation}
  for $t\in I_m$, where
  \begin{equation}
    \begin{split}
      \spaceind[m]^\star(t):=&\
      \cE_{\pivot,\xespace_m}[\full u'(t), (\feellop[m] \full u'
        +\spaceproj[m]\trecfeta{U}''-\trecfeta{U}'')(t)]\\
      &+
      |\varkappa_m (t)|
      \cE_{\pivot,\xespace_m}[\full u(t_{m-1}^+),	(\feellop [m]
        \full u+ \spaceproj[m]\trecfeta{U}'-\trecfeta{U}')(t_{m-1}^+)]
      \\
      &+
      |\varkappa_m (t)|
      \cE_{\pivot,\xespace_{m-1}}[\full u(t_{m-1}^-),(\feellop [m-1] \full u+ \spaceproj[m-1]\trecfeta{U}'-\trecfeta{U}')(t_{m-1}^-)]
      .
    \end{split}
  \end{equation}
  Combining the above estimates, we arrive at 
  \begin{equation}
    \begin{split}
      \Norm{\rho}_{\leb{\infty}(0,t_n;\pivot)}
      &\le
      \sqrt{2}\ltwonorm{u_0-\spaceproj[0] u_0}+2\Norm{\fullproj[]  f-f }_{\leb{1}(0,t_n;\pivot)}\\
      &\quad +2\Big(t_n\sum_{m=1}^n\constref[\timestep m,\tdeg m]{eqn:time-reconstruction:constant}^2
      \Norm{\jumpatof{m-1}{\fullproj[]f-\trecfeta u'}}_{\pivot}^2\Big)^{\nicefrac12}\\
      &\quad +2\sum_{m=1}^n \int_{I_m} \spaceind[m]^\star\d t;
    \end{split}
  \end{equation}    
  this, in conjunction with \eqref{eqn:ell_t_error_linftyl2} now yields
  an alternative \aposteriori error bound which may be superior to the
  one given in Theorem \ref{the:LinfL2final} for small final time $t_n$.
\end{Obs}
\section{$\sobh1(\cI;\dual\elldom)$-type \aposteriori error estimates}
\label{sec:H1(Hm1)apost}
We conclude this work by briefly arguing on how our
techniques allow to derive $\sobh1(\cI;\dual\elldom)$-type
\aposteriori error estimates. To this end, for any $\cI$-piecewise
sufficiently smooth function~$z$, we define the broken (semi-)norm
\begin{equation}
  \tnorm{z}_{\sobh1(\cI;\dual\elldom)}
  :=\Big(
  \sum_{n=1}^N\Norm{\pdt{z}+\chi_n\qp{\jumpatof{n-1}{z}}}^2_{\leb2(I_n;\dual\elldom)}
  \Big)^{\nicefrac12},
\end{equation}
with~$\cI$ signifying the time partition of~$I=(0,T]$ from~\eqref{eq:cI}. Then, recalling~$\rho$ from~\eqref{eq:rho1}, the triangle inequality yields
\begin{equation}
  \label{eq:proof:H1H--1-error-estimate:full-error-split:norm-triangle}
  \tnorm{  u -\full u}_{\sobh1(\cI;\dual\elldom)}
  \leq
  \norm{\rho}_{\sobh1(I;\dual\elldom)}
  +
  \tnorm{w-\full u}_{\sobh1(\cI;\dual\elldom)},
\end{equation}
  upon noting the continuity of $\rho$ with respect to the time
  variable. To control the first term on the right-hand side
of~\eqref{eq:proof:H1H--1-error-estimate:full-error-split:norm-triangle},
we start from~\eqref{eq:spacetime_strong_two}, and notice that
\begin{equation}
  \norm{\rho}_{\sobh1(I;\dual\elldom)}
  =\Norm{\allres{}-\ellop\rho}_{\leb2(I;\dual\elldom)}
  \leq
  \Norm{\allres{}}_{\leb2(I;\dual\elldom)}
  +
  \beta\Norm{\rho}_{\leb2(I;\elldom)},
\end{equation}
with~$\beta$ from~\eqref{eq:ellop-prop}.
Furthermore, based on \eqref{eq:Linfty_mid0} with $t=T$, standard arguments give
  \begin{equation}	\label{eq:Hone_one}
    \alpha \dgenorm[I]{\rho}^2
    \le \ltwonorm{{\rho}(0)}^2+\alpha^{-1}	\Norm{\allres{}}_{\leb2(I;\dual\elldom)}^2.
  \end{equation}   
  Combining the last two bounds, we arrive at 
  \begin{equation}
    \label{eq:proof:H1H--1-error-estimate:full-norm-set-up}
    \norm{\rho}_{\sobh1(\cI;\dual\elldom)}^2 \leq
    \frac{2\beta^2}{\alpha} \ltwonorm{{\rho}(0)}^2 +
    2\Big(1+\frac{\beta^2}{\alpha^{2}}\Big)\Norm{\allres{}}_{\leb2(I;\dual\elldom)}^2.
  \end{equation}
  The first term on the right-hand side
  of~\eqref{eq:proof:H1H--1-error-estimate:full-norm-set-up} is
  bounded trivially by $\eta^{\rm in}$ (defined in Theorem \ref{thm:L2H1}),
  while for the second we use \eqref{eqn:bound_on_xi}.
  To estimate the second term on the right-hand side
  of~\eqref{eq:proof:H1H--1-error-estimate:full-error-split:norm-triangle},
  we consider the splitting
  \begin{equation}
    \tnorm{w - \full U}_{\sobh1(\cI;\dual\elldom)}
    \le
    \tnorm{w - \trecfeta u}_{\sobh1(\cI;\dual\elldom)}
    +\tnorm{ \trecfeta u - \full U}_{\sobh1(\cI;\dual\elldom)}.
  \end{equation}
  The second term on the right-hand side of the last estimate is a
  computable quantity, while the first can be immediately bounded
  using Lemma \ref{lem:ell_rec}, through the linearity of the time
  reconstruction.

\tw{
\section{Conclusions}
In this article we have presented \aposteriori error estimates for fully discrete $hp$-DG-in-time and general Galerkin-in-space discretizations of abstract linear parabolic PDE. Our approach is based on a novel combination of temporal and elliptic reconstructions, the latter allowing for arbitrary elliptic error spatial estimators. Our main results include computable $\leb{\infty}(\pivot)$- and
  $\leb{2}(\elldom)$-\aposteriori error estimates; some remarks concerning $\sobh1(\dual\elldom)$-norm error estimation are given as well. Finally, we note that a series of numerical experiments showcasing the optimality of the \aposteriori error estimators derived above are given in \cite{mohammad}.
}

\subsection*{Acknowledgements}
We would like to thank the referees for their thoughtful remarks that have led
to relevant improvements of our results. EHG acknowledges financial support by The Leverhulme Foundation (Grant No.~RPG-2015-306) and by the Marie Skłodowska--Curie ITN ``ModCompShock''. OL's contribution was partially supported by the Marie Skłodowska--Curie ITN ``ModCompShock''.
TPW acknowledges the financial support of the Swiss National Science Foundation (SNF), Grant No. 200021\underline{\space}182524.

\bibliographystyle{amsalpha}

\begin{thebibliography}{WGSS01}

\bibitem[AMN06]{AMN06}
G.~Akrivis, C.~Makridakis, and R.~H. Nochetto, \emph{A posteriori error
  estimates for the {C}rank-{N}icolson method for parabolic equations}, Math.
  Comp. \textbf{75} (2006), no.~254, 511--531.

\bibitem[AO00]{AinsworthOden:00:book:A-posteriori}
Mark Ainsworth and J.~Tinsley Oden, \emph{A posteriori error estimation in
  finite element analysis}, Pure and Applied Mathematics (New York),
  Wiley-Interscience [John Wiley \& Sons], New York, 2000. \MR{MR1885308
  (2003b:65001)}

\bibitem[BPS09]{BPS:09}
D.~Braess, V.~Pillwein, and J.~Sch\"oberl, \emph{Equilibrated residual error
  estimates are p-robust}, Computer Methods in Applied Mechanics and
  Engineering \textbf{198} (2009), no.~13-14, 1189--1197.

\bibitem[Bra07]{Braess:07:book:Finite}
Dietrich Braess, \emph{Finite elements: theory, fast solvers, and applications
  in elasticity theory}, 3 ed., Cambridge University Press, Cambridge, 04 2007.

\bibitem[BS07]{BrennerScott:07:book:The-mathematical}
S.~C. Brenner and L.~R. Scott, \emph{The mathematical theory of finite element
  methods}, third ed., Texts in Applied Mathematics, vol.~15, Springer-Verlag,
  New York, 2007.

\bibitem[CGM14]{CGM}
A.~Cangiani, E.~H. Georgoulis, and S.~Metcalfe, \emph{Adaptive discontinuous
  {G}alerkin methods for nonstationary convection-diffusion problems}, IMA J.
  Numer. Anal. \textbf{34} (2014), no.~4, 1578--1597. \MR{3269437}

\bibitem[CGS20a]{mohammad}
A.~Cangiani, E.~H. Georgoulis, and M.~Sabawi, \emph{{\it {A} posteriori} error
  analysis for implicit-explicit {$hp$}-discontinuous {G}alerkin timestepping
  methods for semilinear parabolic problems}, J. Sci. Comput. \textbf{82}
  (2020), no.~2, Paper No. 26, 24. \MR{4056776}

\bibitem[CGS20b]{oliver}
Andrea Cangiani, Emmanuil~H. Georgoulis, and Oliver~J. Sutton, \emph{Adaptive
  non-hierarchical {Galerkin} methods for parabolic problems with application
  to moving mesh and virtual element methods}, online preprint 2005.05661,
  arXiv, May 2020.

\bibitem[DM99]{DM:99}
P.~Destuynder and B.~M\'{e}tivet, \emph{Explicit error bounds in a conforming
  finite element method}, Math. Comp. \textbf{68} (1999), no.~228, 1379--1396.
  \MR{1648383}

\bibitem[EJ91]{EJ-I}
K.~Eriksson and C.~Johnson, \emph{Adaptive finite element methods for parabolic
  problems. {I}. {A} linear model problem}, SIAM J. Numer. Anal. \textbf{28}
  (1991), no.~1, 43--77.

\bibitem[EJ95a]{EJ-II}
\bysame, \emph{Adaptive finite element methods for parabolic problems. {II}.
  {O}ptimal error estimates in {$L_\infty L_2$} and {$L_\infty L_\infty$}},
  SIAM J. Numer. Anal. \textbf{32} (1995), no.~3, 706--740.

\bibitem[EJ95b]{EJ-IV}
\bysame, \emph{Adaptive finite element methods for parabolic problems. {IV}.
  {N}onlinear problems}, SIAM J. Numer. Anal. \textbf{32} (1995), no.~6,
  1729--1749.

\bibitem[EJ95c]{EJ-V}
\bysame, \emph{Adaptive finite element methods for parabolic problems. {V}.
  {L}ong-time integration}, SIAM J. Numer. Anal. \textbf{32} (1995), no.~6,
  1750--1763.

\bibitem[EJL98]{EJ-VI}
K.~Eriksson, C.~Johnson, and S.~Larsson, \emph{Adaptive finite element methods
  for parabolic problems. {VI}. {A}nalytic semigroups}, SIAM J. Numer. Anal.
  \textbf{35} (1998), no.~4, 1315--1325 (electronic).

\bibitem[EJT85]{EJT85}
K.~Eriksson, C.~Johnson, and V.~Thomée, \emph{Time discretization of parabolic
  problems by the discontinuous {G}alerkin method}, RAIRO Mod\'el. Math. Anal.
  Num\'er. \textbf{19} (1985), no.~4, 611--643.

\bibitem[ESV17]{ESV:16}
A.~Ern, I.~Smears, and M.~Vohral\'{\i}k, \emph{Guaranteed, locally space-time
  efficient, and polynomial-degree robust a posteriori error estimates for
  high-order discretizations of parabolic problems}, SIAM J. Numer. Anal.
  \textbf{55} (2017), no.~6, 2811--2834. \MR{3723331}

\bibitem[ESV19]{ESV:17}
\bysame, \emph{Equilibrated flux {\it a posteriori} error estimates in
  {$L^2(H^1)$}-norms for high-order discretizations of parabolic problems}, IMA
  J. Numer. Anal. \textbf{39} (2019), no.~3, 1158--1179. \MR{3984054}

\bibitem[GM14]{GM}
E.~H. Georgoulis and C.~Makridakis, \emph{On a posteriori error control for the
  {A}llen-{C}ahn problem}, Math. Methods Appl. Sci. \textbf{37} (2014), no.~2,
  173--179. \MR{3153090}

\bibitem[GSKZ19]{GKSZ17}
F.~D. Gaspoz, K.~Siebert, C.~Kreuzer, and D.~A. Ziegler, \emph{A convergent
  time-space adaptive {${\rm dG}(s)$} finite element method for parabolic
  problems motivated by equal error distribution}, IMA J. Numer. Anal.
  \textbf{39} (2019), no.~2, 650--686. \MR{3941881}

\bibitem[HW18]{HolmWihler:15}
B.~Holm and T.~P. Wihler, \emph{Continuous and discontinuous {G}alerkin time
  stepping methods for nonlinear initial value problems with application to
  finite time blow-up}, Numer. Math. \textbf{138} (2018), no.~3, 767--799.
  \MR{3767700}

\bibitem[Jam78]{J78}
P.~Jamet, \emph{Galerkin-type approximations which are discontinuous in time
  for parabolic equations in a variable domain}, SIAM J. Numer. Anal.
  \textbf{15} (1978), no.~5, 912--928.

\bibitem[KMW18]{KMW:16}
I.~Kyza, S.~Metcalfe, and T.~P. Wihler, \emph{{$hp$}-adaptive {G}alerkin time
  stepping methods for nonlinear initial value problems}, J. Sci. Comput.
  \textbf{75} (2018), no.~1, 111--127. \MR{3770314}

\bibitem[LM06a]{LakMak06}
O.~Lakkis and C.~Makridakis, \emph{Elliptic reconstruction and a posteriori
  error estimates for fully discrete linear parabolic problems}, Math. Comp.
  \textbf{75} (2006), no.~256, 1627--1658.

\bibitem[LM06b]{LakkisMakridakis:06:article:Elliptic}
Omar Lakkis and Charalambos Makridakis, \emph{Elliptic reconstruction and a
  posteriori error estimates for fully discrete linear parabolic problems},
  Math. Comp. \textbf{75} (2006), no.~256, 1627--1658. \MR{2240628
  (2007e:65122)}

\bibitem[MB97]{MB:97}
C.~Makridakis and I.~Babu\v{s}ka, \emph{On the stability of the discontinuous
  {G}alerkin method for the heat equation}, SIAM J. Numer. Anal. \textbf{34}
  (1997), no.~1, 389--401.

\bibitem[MN03]{MakNoc03}
C.~Makridakis and R.~H. Nochetto, \emph{Elliptic reconstruction and a
  posteriori error estimates for parabolic problems}, SIAM J. Numer. Anal.
  \textbf{41} (2003), no.~4, 1585--1594.

\bibitem[MN06a]{MakNoc06}
\bysame, \emph{A posteriori error analysis for higher order dissipative methods
  for evolution problems}, Numer. Math. \textbf{104} (2006), no.~4, 489--514.

\bibitem[MN06b]{MakridakisNochetto:06:article:A-posteriori}
Charalambos Makridakis and Ricardo~H. Nochetto, \emph{A posteriori error
  analysis for higher order dissipative methods for evolution problems}, Numer.
  Math. \textbf{104} (2006), no.~4, 489--514. \MR{2249675 (2008b:65114)}

\bibitem[MSW05]{MSW1}
A.-M. Matache, C.~Schwab, and T.~P. Wihler, \emph{Fast numerical solution of
  parabolic integrodifferential equations with applications in finance}, SIAM
  J. Sci. Comput. \textbf{27} (2005), 369--393.

\bibitem[MSW06]{MSW2}
\bysame, \emph{Linear complexity solution of parabolic integro-differential
  equations}, Numer. Math. \textbf{104} (2006), 69--102.

\bibitem[Pic98]{P98}
M.~Picasso, \emph{Adaptive finite elements for a linear parabolic problem},
  Comput. Methods Appl. Mech. Engrg. \textbf{167} (1998), no.~3-4, 223--237.

\bibitem[Rou13]{Roubicek:13}
T.~Roub{\'{\i}}{\v{c}}ek, \emph{Nonlinear partial differential equations with
  applications}, second ed., International Series of Numerical Mathematics,
  vol. 153, Birkh\"auser/Springer Basel AG, Basel, 2013.

\bibitem[SS00]{SchoetzauSchwab00}
D.~Schötzau and C.~Schwab, \emph{Time discretization of parabolic problems by
  the $hp$-version of the discontinuous {G}alerkin finite element method}, SIAM
  J. Numer. Anal. \textbf{38} (2000), 837--875.

\bibitem[SS01]{SchoetzauSchwab01}
\bysame, \emph{$hp$-discontinuous {G}alerkin time-stepping for parabolic
  problems}, C. R. Acad. Sci. Paris, S\'erie I \textbf{333} (2001), 1121--1126.

\bibitem[SW10]{SW10}
D.~Sch\"otzau and T.~P. Wihler, \emph{A posteriori error estimation for
  {$hp$}-version time-stepping methods for parabolic partial differential
  equations}, Numer. Math. \textbf{115} (2010), no.~3, 475--509.

\bibitem[Ver98]{V98}
R.~Verfürth, \emph{A posteriori error estimates for nonlinear problems.
  {$L^r(0,T;L^\rho(\Omega))$}-error estimates for finite element
  discretizations of parabolic equations}, Math. Comp. \textbf{67} (1998),
  no.~224, 1335--1360.

\bibitem[Ver03]{V03}
\bysame, \emph{A posteriori error estimates for finite element discretizations
  of the heat equation}, Calcolo \textbf{40} (2003), no.~3, 195--212.

\bibitem[vPS04]{SchwabPetersdorff}
T.~von Petersdorff and C.~Schwab, \emph{Numerical solution of parabolic
  equations in high dimensions}, Math. Model. Anal. Numer. \textbf{38} (2004),
  93--127.

\bibitem[WGSS01]{Gerdes}
T.~Werder, K.~Gerdes, D.~Sch\"otzau, and C.~Schwab, \emph{$hp$-discontinuous
  {G}alerkin time-stepping for parabolic problems}, Comput. Methods Appl. Mech.
  Engrg. \textbf{190} (2001), 6685--6708.

\end{thebibliography}
\providecommand{\bysame}{\leavevmode\hbox to3em{\hrulefill}\thinspace}
\providecommand{\MR}{\relax\ifhmode\unskip\space\fi MR }
\providecommand{\MRhref}[2]{%
  \href{http://www.ams.org/mathscinet-getitem?mr=#1}{#2}
}
\providecommand{\href}[2]{#2}

\end{document}